\documentclass{svjour3}                     
\smartqed  
\usepackage{graphicx}
\usepackage[round]{natbib}

\usepackage{hyperref}
\hypersetup{colorlinks=true,citecolor=blue,urlcolor=red}
\usepackage[english]{babel}

\usepackage{amssymb}
\usepackage{amsmath}
\usepackage{amsfonts}

\newcommand{\blambda}{{\boldsymbol \lambda}}
\def\by{{\boldsymbol y}}
\def\bh{{\boldsymbol h}}
\def\bA{{\boldsymbol A}}
\def\bB{{\boldsymbol B}}
\def\bT{{\boldsymbol T}}
\def\bx{{\boldsymbol x}}
\def\b1{{\boldsymbol 1}}
\def\bX{{\boldsymbol X}}
\def\bu{{\boldsymbol u}}
\def\bw{{\boldsymbol w}}
\def\bS{{\boldsymbol S}}
\def\bL{{\boldsymbol L}}
\def\bxi{{\boldsymbol \xi}}
\def\bphi{{\boldsymbol \phi}}

 \journalname{Computational Economics}
\begin{document}
\bibliographystyle{plainnat}

\title{Portfolio optimization  in incomplete markets and price constraints determined by maximum entropy in the mean}

\titlerunning{Portfolio optimization by maximum entropy in the mean}

\author{Argimiro Arratia          \and
        Henryk Gzyl
}

\authorrunning{A. Arratia, H. Gzyl}

\institute{A. Arratia \at
              Computer Science, Barcelona Tech (UPC), Barcelona, Spain\\
                           \email{argimiro@cs.upc.edu},\ ORCID:  0000-0003-1551-420X         
           \and
	H. Gzyl \at
	Center for Finance, IESA, Caracas, Venezuela\\
 \email{henryk.gzyl@iesa.edu.ve}
 }

\date{Cite this article:
Arratia, A., Gzyl, H. Portfolio Optimization in Incomplete Markets and Price Constraints Determined by Maximum Entropy in the Mean. {\em Computational Economics},  56, 929-952 (2020). https://doi.org/10.1007/s10614-019-09954-3
}

\maketitle

\begin{abstract}
 A solution to a portfolio optimization problem is always conditioned by constraints on the initial capital  and the price of the available market assets.  If a risk neutral measure is known, then the price of each asset is the discounted expected value of the asset's price under this measure. 
 But if the market is incomplete, the risk neutral measure is not unique, 
 and there is a range of possible prices for each asset, which can be 
 identified with bid-ask ranges.
 We present in this paper an effective method to determine the cu\-rrent 
 prices of a collection of assets  in incomplete markets,
  and such that these prices  comply with the cost constraints for a portfolio optimization problem. 
Our workhorse is the method of maximum entropy in the mean to adjust a distortion function from bid-ask market data.
This distortion function plays the role of a risk neutral measure, which is used to price the assets, and the distorted probability that it determines reproduces bid-ask market values. We carry out numerical examples to study the effect on portfolio returns of the computation of prices of the assets conforming the portfolio with the proposed methodology. 

\end{abstract}

\section{Introduction}
\label{sec:intro}
There are many ways of optimizing portfolios. One can safely say that all of them, except for the naive portfolio diversification, consist of some form of trade off between return and risk subject to an available capital constraint. One of the  standard versions of the portfolio optimization problem, originally due to  \cite{Mark52}, has many variants, all of which use the volatility (or the standard deviation) of the portfolio as risk measure, and the trade off between risk and return is set up in a direct form, like for example: find the portfolio of minimum variance that yields  a given expected return. This last constraint is clearly equivalent to a cost constraint upon the future values of the assets when a risk neutral probability is available.  A few standard references for this matter are  \cite{Plis,KS,FS}. 

Motivated by the impossibility of using non-arbitrage arguments when the market is incomplete, or when the class of risk neutral measures is not known, or more importantly, when the initial price for each asset can not be specified, an alternative theoretical methodology, which has been termed ``conic portfolio theory'', was developed recently  by \cite{Madan16} (see also \cite{MS}),  and others. In this proposal, ideas coming  from risk analysis in 
 the insurance industry are used to provide the trade off  between return and risk, by letting  the role of the risk constraint to be played by either the bid or ask prices of the assets, which are estimated using a distortion function.

Here we shall use distortion functions in a different, almost opposite way. For this we think of the price of each asset as the price of a collection of risks specified up to a bid-ask range, and then we determine a risk pricing measure that yields prices within the bid-ask ranges. The prices of assets  so determined will be the current  prices which are to be used to compute the returns of the portfolio.

The  numerical procedure we propose to obtain the  distortion function is based on the method of maximum entropy in the mean. This is a convenient method to solve inverse problems subject to convex constraints upon the solution, which allows for the data to be specified in intervals. This procedure has the important advantage of being model free, which in particular does not call for the calibration of parameters.  

The historical antecedents of the maximum entropy in the mean method can be traced back to 
\cite{jaynes1957}, who proposed and used a standard  maximum entropy  methodology in statistical physics. This form of the maximum entropy method is later used as a valuable approximation scheme to solutions of some cases of the problem of moments in \cite{mead1984}, and the mathematical foundation of the maximum entropy in the mean method is developed in \cite{dacunha1990}. 

Here the maximum entropy in the mean method is used to assign a price to a market asset, which in a certain sense corresponds to the discounted expected value of the asset under some (unknown) risk-neutral measure, but what can only be said  is that this price is within some market  bid-ask range.
We termed these prices, determined by the market, as {\em conservative market prices}, in concordance with similar terminology use  in the literature \citep{Madan16,MS}.
Once these conservative market  prices for the assets comprising our portfolio are at hand, we can do two things. First, we can use the  prices to compute returns, and second, to assign value to the price constraint for the portfolio optimization problem. 

Interrelationships between derivative valuation in incomplete markets and risk pricing in insurance have been brought up many times. When the market model is incomplete, it was noted that instead of a given price, one could assign a price range to any given asset. One could also use methods from utility (or risk) theory to determine a market price of the asset. 
In the insurance industry a similar problem exists: How to assign a price to a potential loss, or a price to risk. For this purpose, a risk pricing axiomatics and a risk pricing methodology were developed, which in many ways parallels and overlaps that developed to value the risk of financial assets. For a few references about this see  \cite{Y} and  \cite{LG}.

\subsection{Plan of our work}
We begin with briefly describing the market models in finance and in insurance,
at the same time that we introduce the necessary notation in Section \ref{basics}. 
There are two issues related to bid-ask ranges. On one hand, bid-ask ranges for assigning prices to derivatives are related to market incompleteness. When the market is incomplete but the class of risk neutral measures is known, a bid-ask range can be applied
 (but see \cite{CH00} 
 for a characterization in term of no good deal prices).
On the other hand, the existence of bid-ask ranges for current prices does not allow the construction of the class of risk neutral probabilities, which impedes risk neutral pricing. But equally important, since to solve the problem of portfolio optimization we need to be able to calculate asset returns, we have to device a way of going around the absence of current asset price.

We briefly describe the way this is taken care of in \cite{MS}, 
\cite{Madan16}. Since their proposal involves the use of distorted probabilities and their connection to coherent risk measures, we devote some space to coherent risk measures, risk pricing in insurance, and eventually establish a connection between distorted risk measures and risk neutral probabilities, that can be summarised as follows (see Section \ref{DF_RNP}). 
We take the bid-ask range as the only available market data, and determine a distortion function that yields asset prices within the bid-ask range. 
Since the distorted probability law so obtained is equivalent to a risk neutral probability yielding the computed price as actual price, we shall consider that price as the base price against which we compute the returns of the basic assets for the purpose of portfolio optimization.

In Section \ref{discret} we describe the numerical procedure to obtain the distortion function, show how the discretization is carried out, and write down  the solution of the problem by the method of maximum entropy in the mean. It will turn out that there will be two possible solutions, depending on whether we do, or do not, impose bounds upon the derivatives of the distortion function. We present both results there, but we work out the details of this procedure in the Appendix \ref{A:sec1}. 
 
We consider in  Section \ref{portfolio} a portfolio optimization problem. For that we shall consider an exponential utility and make use of the connection between that problem and entropy maximization, in a way similar to that developed by  \cite{FS}.

In  Section  \ref{conc} we sum it all up, and  in the Appendices \ref{A:sec1} and \ref{A:sec2} we collect all material set aside for not to interrupt the main flow of ideas.
  
\section{Bid-ask  prices and distortion functions}\label{basics}
Before entering into the subject of the section, it is convenient to establish the market models in finance and insurance in order to introduce the basic notations.
\subsection{The market models}
In this work we shall consider a two time market model consisting of a probability space $(\Omega,\mathcal{F},P)$ and a finite collection of $M$ basic assets whose prices at time $t=0$ are denoted by 
$S_1(0), \ldots, S_M(0)$ and whose prices at time $t=1$ are positive, continuous random variables 
$S_1(1), \ldots, S_M(1)$,
each distributed according to a probability density function $F_i$, $i = 1, \ldots, M$.

 In this model we might as well suppose that $\Omega=[0,\infty)^{M},$ and that $\mathcal{F}=\mathcal{B}(\bS(1)),$ the Borel $\sigma$-algebra generated by 
$\bS(1)=(S_1(1), \ldots, S_M(1)).$
Here $P$ denotes a probability measure on $(\Omega, \mathcal{F}),$ the so called ``physical law'' which describes the statistical behavior of $\bS(1)$. 
 A standard financial asset in this model will be a finite valued, 
 $\mathcal{F}$-measurable function  $X:\Omega\rightarrow\mathbb{R},$ and the problem for the financial analyst is how to price $X$.

The simplest model used in the insurance industry goes as follows. We start with a probability space 
$(\Omega,\mathcal{F},P)$,  and a collection of positive random variables denoting the severity  (or monetary value) of the losses occurring in a given period of time. We do not have the analogue of the initial prices, just the loss incurred in the period. And we suppose that there is a finite number $L_1, \ldots,L_M$ of them, denoted collectively by $\bL.$ As before, we suppose that $\Omega=[0,\infty)^M,$ and that $\mathcal{F}=\mathcal{B}(\bL)$. In this case we are interested in pricing the losses or functions of the losses, which can be thought as the prices of the possible insurance products. 
 
\subsection{Bid-ask prices}
A key tool for the pricing of financial assets is that of risk neutral probability. The class of such probabilities is defined (in a one period) market model by
$$\mathcal{M}_e = \{Q\in\mathcal{M}(P)\,|\,E_Q[e^{-r}S_i(1)] = S_i(0)\}\;\;\;i=1, \ldots, M.$$
A market model is complete when this class reduces to one element. Otherwise, the market model is said to be incomplete, and financial assets do not have a unique price. The theoretical determination of the bid-ask price of an asset of future value $X$ at time $t=1$ is given by (see \cite{Plis} or  \cite{KS}):
\begin{eqnarray}\label{ba}
b(X) &=& \exp(-r) \inf \{E_Q[X] : Q \in {\cal M}_e \} \nonumber\\
a(X) &=& \exp(-r) \sup \{E_Q[X] : Q \in {\cal M}_e \}
\end{eqnarray}
As mentioned, the class $\mathcal{M}_e$ may not be known, then in this case (\ref{ba}) can not be used for obtaining a price range. One reason for which $\mathcal{M}_e$ may not be known is that the basic assets do not have a specified price at $t=0,$ but are known only up to a bid-ask range. That is, we are in what is called a two price economy in \cite{Madan16} or \cite{MS}. These authors propose a different market model in which there are no basic assets, and no notion of risk neutral measure.

To be specific, they propose a notion of bid-ask range that goes as follows. To begin with consider a convex cone $\mathcal{A}$ of acceptable risks, and a class of probability laws $\mathcal{P}=\{Q \sim P| E_Q[X] \geq 0;\;\forall X\in\mathcal{A}\}$ (see Chapter IV of \cite{FS} for these notions), and propose the following definition of bid-ask range for each $X$ as follows:
\begin{eqnarray}\label{ba2}
b(X) &=&  \inf \{E_Q[X] : Q \in \mathcal{P} \} \nonumber\\
a(X) &=&  \sup \{E_Q[X] : Q \in \mathcal{P}\}
\end{eqnarray}
 
 At this point we mention again that here we shall regard the bid-ask range of the basic assets as a datum obtained from the market, either as quoted quantities or to be determined from the daily price data. 
 
\subsection{Coherent and convex risk measures}
The need to appropriately price risk resulted in an rather precise (axiomatic) characterization of what is a risk measure. Actually, there already exist several interrelated classes of risk measures, which are already a text book material, and what follows about coherent and convex risk measures is taken from   \cite{FS}, 
although originally proposed by  \cite{artzner1999}.
Let $\Xi$ be the collection of random variables representing the changes of value of an asset. A {\em risk measure} is any function $\rho:\Xi\rightarrow \mathbb{R},$ and a {\em coherent risk measure} is defined as follows.
\begin{definition}\label{cohrm}
A coherent risk measure is any risk measure satisfying the following conditions:
\begin{itemize}
\item[{\bf 1.}] (Monotonicity) For $X_1,X_2 \in \Xi,$ and $X_1\leq X_2$ then $\rho(X_2) \leq \rho(X_1).$
\item[{\bf 2.}] (Sub-additivity) For $X_1,X_2 \in \Xi,$ it holds that $\rho(X_1 + X_2) \leq \rho(X_1)+ \rho(X_2).$
\item[{\bf 3.}] (Homogeneity) For any $a \geq 0$ and any $X\in\Xi,$ we have $\rho(aX) = a\rho(X).$
\item[{\bf 4.}] (Translation invariance) For any $m\in\mathbb{R}$ and any $X\in \Xi,$ $\rho(X+m)=\rho(X)-m.$
\end{itemize}
\end{definition}
When the sub-additivity and homogeneity conditions are merged into the convexity requirement:
$$\mbox{For}\;\; X_1,X_2 \in \Xi\;\;\mbox{and}\;\; t\in [0,1],\;\; \rho(tX_1 +(1-t)X_2) \leq t\rho(X_1) + (1-t)\rho(X_2)$$
\noindent we say that $\rho$ is a {\em convex risk measure}.

The following result, which can be derived from \cite[Prop. 4.1]{artzner1999},
 suggests a similarity between bid-ask prices and coherent risk measures.
\begin{theorem}\label{dualrep}
A lower semi-continuous risk measure $\rho$ is coherent if and only if there exists   a class $\mathcal{Q}\subset\mathcal{M}(P)$ such that
$$\rho(X) = \inf\{E_Q[X]\,|\,Q \in \mathcal{Q}\}.$$
\end{theorem}

The probabilities in $\mathcal{Q}$ are called ``generalized scenarios''.
Observe that, according to (\ref{ba}), when the class of generalized scenarios is the class of risk neutral measures, this theorem renders the bid  price $b(X)$ as a risk measure  (if $X$ could be thought of as a change of value of some asset). Similarly for $-a(-X).$  
This analogy is behind the portfolio optimization methodology propounded in  \cite{Madan16}.

\subsection{Risk pricing measures}
Risk pricing in insurance exists from a long time ago, and many ways to compute premiums exist. 
Eventually various premium principles were established, and depending on the premium principles chosen,  
 computational methods consistent with them were proposed. See   \cite{WYP},   \cite{Y} or   \cite{LG}, for example. To stress the analogies with risk axiomatics, let us list the premium principles.
\begin{definition}\label{preprin}
Let $\Xi$ be a collection of positive, integrable random variables modeling losses. A premium principle is a mapping $\Pi:\Xi\rightarrow [0,\infty]$ satisfying:
\begin{itemize}
\item[{\bf 1.}] (Law invariance) $\Pi(X)$ depends only on $F_X.$
\item[{\bf 2.}] (Homogeneity) For any $X\in \Xi$ and $a>0,$ we have $\Pi(aX) = a\Pi(X).$
\item[{\bf 3.}] (Translation invariance) For any $X\in \Xi$ and $a>0,$ we have $\Pi(a + X) = a + \Pi(X).$\\
Conditions 2 and 3 give the ({\it no unjustified risk loading}) principle that states that the price of a known loss is the amount of the loss; $\Pi(a)=a$ for any $a\geq 0.$
\item[{\bf 4.}] (Risk loading) For any $X\in \Xi,$ $\Pi(X)\geq E[X].$
\item[{\bf 5.}] (No rip-off) For any $X\in \Xi,$ $\Pi(X)\geq esssup(X).$
\item[{\bf 6.}] (Sub-additivity) For any $X_1, X_2\in \Xi$ it holds that $\Pi(X_1 + X_2) \leq \Pi(X_1)+\Pi(X_2).$ Additivity and super-additivity can be similarly defined.
\item[{\bf 7.}] (Additivity for comonotonic losses). Two random losses $X_1, X_2\in \Xi$ are comonotonic if for every $\omega_1, \omega_2 \in \Omega$ we have $(X_1(\omega_2)-X_1(\omega_1))(X_2(\omega_2)-X_2(\omega_1))\geq 0.$ For comonotonic variables, $\Pi$ is comonotonic additive if $\Pi(X_1 + X_2) = \Pi(X_1) + \Pi(X_2).$
\item[{\bf 8.}] (Monotonicity) For $X_1, X_2\in \Xi$ such that $X_1 \leq X_2,$ then $\Pi(X_1) \leq \Pi(X_2).$
\item[{\bf 9.}] (Stochastic ordering compatibility) A premium is said to be compatible with first sto\-chas\-tic dominance if $F_{X_1} \leq F_{X_2},$ then $\Pi(X_1) \leq \Pi(X_2).$ Similarly, it is said to be second order stochastic dominance compatible, if whenever $E[(X_1-d)^+]\leq E[(X_2-d)^+]$, for all $d>0,$ then $\Pi(X_1) \leq \Pi(X_2).$\\
(Here $(X-a)^+$ denotes as usual $\max(X-a,0)$.)
\end{itemize}
\end{definition}

The following result from   \cite{WYP} 
ties coherent risk theory with risk pricing in insurance and will be the basis of our proposal to obtain an alternative to risk neutral measures from bid-ask prices. 
(Recall that $X \wedge a$ stands for $\min(X,a)$.)
\begin{theorem}\label{wyp}
Suppose that a risk pricing measure satisfies Properties (1), (2), (3),(7), (8) and the following continuity property:
$$\mbox{For } X\in \Xi,\;\;\;\lim_{a\downarrow 0}\Pi((X-a)^+)=\Pi(X),\;\;\;\lim_{a\downarrow 0}\Pi(X\wedge a)=\Pi(X).$$
Then there exists a non-decreasing (distortion) function $g:(0,1)\rightarrow(0,1)$ such that
\begin{equation}\label{dist1}
\Pi(X) = \int_0^\infty g(\widehat{F}_X(x))dx.
\end{equation}
Here $\widehat{F}_X=1-F_X.$ Also, if $\Xi$ contains all Bernoulli random variables, $g$ is unique and $g(p)$ is given by the price of the Bernoulli(p)  random loss (remember law invariance).
\end{theorem}
Another proof of this representation result appears in   \cite{WW}.

\subsection{Relation to  bid and ask prices}
If we think of the prices $\{S_i(1), i=1, \ldots,M\}$ of the basic assets at time $t=1$ as positively valued risks, and think of the bid-ask price ranges $(b_i,a_i)$ as being given by the market, we may ask ourselves: What is the distortion function producing these prices?
 That is, we want to solve the following problem:  
 \begin{quote}
 Determine a continuously differentiable concave distortion function $g(u)$ such that
\begin{equation}\label{prob0}
e^{-r}\int_0^\infty g(\widehat{F}_{S_i}(x))dx \in (b_i,a_i) , \;\; i=1, \ldots,M.
\end{equation}
\end{quote}
Hence forward  to simplify the notation we use $F_i$ instead of $F_{S_i}.$ Let $\psi(u) =1-g(1-u)$. 
Then 
$$\int_0^\infty xd\psi(F_i(x)) = \int_0^\infty\left(\int_0^\infty I_{[0,x]}(s)ds\right)d\psi(F_i(x))$$ 
$$ = \int_0^\infty\Big(\int_0^\infty I_{[s,\infty}(x)d\psi(F_i(x))\Big)ds = \int_0^\infty g(\widehat{F}_i(s))ds.$$
But also in account of the assumptions on the $S_i$ and $g,$ making the change of variables $u=F_i(x)$ or $x=F_i^{-1}(u)=q_i(u),$ we obtain
\begin{eqnarray*}
\int_0^\infty xd\psi(F_i(x)) &=& \int_0^\infty xg'(\widehat{F}_i(x))d\widehat{F}_i(x) \\
&=& \int_0^1q_i(u)g'(1-u)du = \int_0^1q_i(u)\phi(u)du.
\end{eqnarray*}
Here we set $\phi(u)=g'(1-u).$ Since we want $g$ to be concave, we have to require $\phi(u)$ to be increasing. Thus with these changes of notation problem (\ref{prob0}) becomes:
\begin{quote}
Determine a continuous, increasing $\phi(u)$ on $[0,1]$ such that
\begin{equation}\label{prob1}
e^{-r}\int_0^1 q_i(u)\phi(u)du \in (b_i,a_i), \;\; i=1, \ldots, M.
\end{equation}
\end{quote}

\subsection{Distortion functions and risk neutral densities}\label{DF_RNP}
Here we establish a relation between the two valuation methods. Let us consider the case of a market model with one random asset and denote by $F$ the cumulative distribution function of its future value $S(1)$. Let $h(x)$ be either a positive or a bounded measurable function defined on $[0,\infty).$ 
 We define the probability $P^{\psi}$ on $(\Omega, \sigma(S(1)))$ by
$$E^{\psi}[h(S(1))] = \int_0^\infty h(x)d\psi(F(x)) = \int_0^\infty h(x)\psi'(F(x))dF(x).$$
Let us put $\Lambda=\psi'(F(S(1)))$.  Observe that  when $\psi$ is concave and increasing $\psi>0,$ thus $\Lambda >0.$ Now define $Q\sim P$ by $dQ/dP = \Lambda$.
  Then clearly
$$E_Q[h(S(1))]=E[\Lambda h(S(1))] = \int_0^\infty h(x)\psi'(F(x))dF(x) = E^{\psi}[h(S(1))].$$
To obtain a distortion function from a risk neutral density, we do as follows. 
Suppose $Q\sim P$ is a risk neutral probability. Let $k$ be a positive Borel measurable function such that $k(S(1))=E[dQ/dP|S(1)].$ Notice that $k(x)$ can be obtained as the Radon-Nikodym derivative of $Q(S(1)\leq x)$ with respect to $P(S(1)\leq x)$ or as the regular conditional probability  $E[dQ/dP|S(1)=x].$ Now define $\psi_1:[0,1]\rightarrow [0,\infty)$ by $\psi_1(u)=k(q(u))$ where $q=F^{-1}$, and
 set $\psi(t)=\int_0^t\psi_1(u)du$.
  Then $\psi'=\psi_1$ for almost all $u\in[0,1]$. 
Notice that $\psi'(F(S(1)))=k(S(1))$ almost surely with respect to the Lebesgue measure on $[0,\infty]$ due to the assumptions on $S(1).$ Furthermore, as above
       $$E^{\psi}[h(S(1))] = E_Q[h(S(1)].$$

Observe that from the interpretation of $k$ mentioned above, it is clear that $\psi(1)=1.$ Now set $g(u)=1-\psi(1-u).$ With the notations as just introduced, we can gather these comments under the statement of the following theorem.
\begin{theorem}\label{RNP-DF}
Let $g(u)$ be a concave, twice continuously differentiable distortion function satisfying 
\[ e^{-r}\int_0^\infty g(\widehat{F}(x))dx \in (b,a) \]
with $b$ and $a$ bid and ask prices at time $t=0$.
 Let $\psi(u)=1-g(1-u)$, and
define $\Lambda=\psi'(F(S))$ so that $ dQ=\Lambda dP$.
  Then $Q\sim P$ and
$$e^{-r}E_Q[S(1)] \in (b,a).$$
Furthermore, for $h$ as above,
$$E^\psi[h(S(1))] = \int_0^\infty h(x)d\psi(F(x))=\int_0^\infty h(x)\psi'(F(x))dF(x)=E_Q[h(S(1))].$$
Conversely, if $Q\sim P$ satisfies $e^{-r}E_Q[S(1)] \in (b,a),$ define $k(x)=dQ(S(1)\leq x)/dF(x),$ and $\psi(t)=\int_0^tk(q(u))du.$ Since $\psi'(F(S(1))=k(S(1)),$ clearly
$$E^\psi[h(S(1))] = E_Q[h(S(1))]$$
\noindent
and furthermore
$$e^{-r}E^\psi[h(S(1))] \in (b,a).$$
\end{theorem}
The proof of this theorem is given in the Appendix \ref{A:sec2}.

\section{Discretization of problem (\ref{prob1})}\label{discret}
Observe that either of problems (\ref{prob0}) or ({\ref{prob1}) is an integral equation with convex constraints upon the solution and with finite data given in an interval. In order to solve them numerically, the first step is to discretize the problem. Let us show how to discretize (\ref{prob1}). Fist we partition $[0,1]$ into $N$  intervals of equal length  
and put $\phi_j=\frac{1}{2}\big(\phi((j-1)/N)+\phi(j/N)\big),$ for $j=1, \ldots, N.$
Now put $B(i,j) = \frac{e^{-r}}{N}q_i(j/N)$, for $i=1, \ldots,M,$ and $j=1, \ldots, N.$ 
Then, problem (\ref{prob1}) becomes: 

\begin{eqnarray}\label{prob2}
\mbox{Solve the system}\qquad \bB\bphi \in \mathcal{K}=\prod_{i=1}^M[b_i,a_i]\\\nonumber
\mbox{Subject to}\qquad 0 < \phi_1 < \phi_2 <  \ldots <\phi_N.\nonumber
\end{eqnarray}
To change the nature of this constraint we proceed as follows. Let $\bT$ be the $M \times N$ lower triangular matrix with all entries equal to $1$ and let $\bxi\in [0,\infty)^N,$  put $\bphi=\bT\bxi$ and $\bA=\bB\bT.$ Then (\ref{prob2}) becomes
\begin{eqnarray}\label{prob3}
\mbox{Solve the system}\qquad \bA\bxi \in \mathcal{K}=\prod_{i=1}^M[b_i,a_i]\\\nonumber
\mbox{Subject to}\qquad \bxi \in [0,\infty)^N.\nonumber
\end{eqnarray}
This is an under-determined linear system with convex constraints on the solution and data in intervals. The method of Maximum Entropy in the Mean (MEM for short) takes care nicely of this type of problems.

\subsection{Estimating the distortion functions by Maximum Entropy in the Mean}
Actually, MEM allows for flexibility in the choice of the constraints. We could consider for example, that either  $\bxi\in[0,L]^N,$ or even in some more elaborate domain, or as above, $\bxi\in[0,\infty).$ The two choices lead to somewhat different representation of the solutions.

Here we shall only list the results and postpone the details to the appendices section. When we require the solution to problem (\ref{prob1}) to be bounded (we shall refer to this case as the MEM bounded), the MEM will provide us with a representation like
\begin{equation}\label{rep01}  
\xi_j^* = L\frac{e^{-L(\bA^t\blambda)_j}}{1 + e^{-L(\bA^t\blambda)_j}},\;;\;j=1,....,N.
\end{equation}
Here $\bA^t$ stands for the transpose of the matrix $\bA,$ and $\blambda^*$ is a particular value of a Lagrange multiplier that we shall specify how to determine in the appendices.
Once $\bxi^*$ is obtained by any method, we can determine the value $\bS(0)$ that enters in (\ref{distpri}) according to
\begin{equation}\label{distpri0}
 \bS(0) = \bA\bxi^* \in \prod_{i=1}^M[b_i,a_i].
\end{equation}

\section{Portfolio optimization}\label{portfolio}
The approach to the portfolio optimization  problem that we are considering consists on constructing a portfolio which maximizes some utility function subject to some initial capital constraints. The utility function is supposed to incorporate the balance between return and risk of the investor.

Therefore, given a utility function $U,$  the physical measure $P$  and the  risk neutral measure $Q$,
 a standard presentation of the portfolio optimization problem  consists in solving the following 
optimization problem (cf.   \cite{FS}): 
 \begin{eqnarray}\label{PO}
 &\mbox{Find}: & \sup_{\bh} \ E^P[U(\langle \bh,\bS(1)\rangle)]  \nonumber\\
 & \mbox{Subject to}: & \langle \bh,\bS(0)\rangle = C_0;\;\;\; h_i \ge 0,  \;\;i=1, \ldots,M.
\end{eqnarray}
where we are considering $M$ assets, whose prices at time $t=0$ are 
$\bS(0) = e^{-r}E_Q[\bS(1)],$  where as above, $\bS(1)$ denotes the random prices at $t=1,$ and where
$C_0$ is the available capital and, for the sake of making the analysis simpler (but without loss of generality),  we  consider only portfolios $\bh = (h_1, \ldots, h_M)$ with long positions in the assets.
 
Alternatively, if one wants to work with unit-less weights, we may consider the gross return of the assets, defined by $X_i(1) = S_i(1)/S_{i}(0), i=1, \ldots, M,$ and setting $w_i = h_iS_i(0)/C_0,$ we can rewrite the portfolio optimization problem as
 \begin{eqnarray}\label{POnounit}
 &\mbox{Find}: & \sup_{\bw} \ E^P[U(C_0\langle \bw,\bX(1)\rangle)]  \nonumber\\
 & \mbox{Subject to}: & \langle \bw,\b1\rangle = 1;\;\;\;0 \leq w_i, \;\;i=1, \ldots,M.
\end{eqnarray}  

Observe that this setup presupposes that we know the prices $S_i(0)$ of the random assets at $t=0.$ Not only that, the sheer notions of gross return or rate of return, imply that we know an initial price of each asset. 

To overcome these two issues we propose to use a distortion function $\Psi$ to define 
{\em conservative market prices} $S^d_i(0)$  by
\begin{equation}\label{distpri}
 S_i^d(0) = e^{-r}\int x d\Psi (F_{S_i}(x))  \in (b_i(x),a_i(x)),\ i=1,\ldots,M
\end{equation}
\noindent as specified in (\ref{distpri0}), and use them to restate the portfolio optimization problem as follows:

 \begin{eqnarray}
 &\mbox{Find}:& \sup_{\bh} \ E^P[U(\langle\bh,\bS(1)\rangle)]  \nonumber\\
 & \mbox{Subject to}:&  \langle \bh,\bS^d(0)\rangle  = C_0;\;\;h_i\geq 0, \;i=1, \ldots,M.\\\nonumber
\end{eqnarray}

As explained in the previous section, we do not need to know (or impose) a particular distortion function, but instead we estimate this from bid and ask prices by the method of maximum entropy in the mean.

\subsection{An  example of portfolio optimization}
In this section we carry out the portfolio optimization procedure for a particular example of utility function. This function is chosen because some of the details of the process overlap those of the entropy maximization procedure for the determination of the distortion function. 
We shall consider the utility function $U(c) = 1 - \exp(-c)$, and then 
to find
\[ \sup_{\bh} \ E^P[U(\langle\bh,\bS(1)\rangle)] \] 
which is equivalent to find 
\[ \inf_{\bh} E^P[ \exp(-\langle\bh,\bS(1)\rangle)]\]
subject to $\langle\bh,\bS^d(0)\rangle  = C_0$, 
where (to insist) $S_i^d(0) = e^{-r}\int x d\Psi (F_{S_i}(x)), \ i=1,\ldots, M$, and we know that $S_i^d(0)\in (b_i,a_i).$ 
The function to minimize is log-convex, so it is rather convenient to write  
$V(\bh) =  E^P[ \exp(-\langle\bh,\bS(1)\rangle)]$ and consider the following equivalent convex minimization problem:

 \begin{eqnarray}\label{utmin1}
 &\mbox{Find}:& \inf_{\bh} \ln V(\bh)  \nonumber\\
  & \mbox{Subject to}:& \langle \bh,\bS^d(0)\rangle = C_0 ;\;\;h_i\geq 0, \;i=1, \ldots,M.
\end{eqnarray}

\subsection{Numerical Experiments}
In order to show the merits of our proposed portfolio optimization model (\ref{utmin1}) in practice, we build a portfolio with real  data for seven companies trading in  the NYSE market, chosen randomly and subject only to have a sufficiently long  daily price history for making extensive simulations. 
The  tickers of the  companies considered are ALB, AAPL, APC, GE, WFC, HAL, and SHW, and for each we have daily OHLC prices history from 2000-01-01 to 2016-12-29.

For our experiments we fix the initial capital $C_0 = 10^5$ units of the market currency.  
 We sampled monthly returns and observed prices    
 at the end of a 
  month, for each asset, 
 so that $\bS(1)$ represents the assets' prices after one month, 
and our goal is to compare the capital gain (or gross return) of the portfolio 
$\bh =(h_1, \ldots, h_M)$  at $t=1$, as given by
$\displaystyle  \frac{\langle \bh,\bS(1)\rangle}{C_0}$, when $\bh$ is obtained from the 
MEM (bounded or unbounded) estimation  $\bS^d(0)$ of the initial prices (i.e. as solution of Equation \ref{utmin1}), 
against the portfolio $\bh$ obtained from current initial price, denoted by  $\widehat{\bS}(0)$. 
We stress the fact that we are doing single period optimization at different epochs to assess our portfolio optimization method in different market scenarios. We are not doing multi-period rebalancing.

\subsubsection{Estimating the bid and ask prices from  data} 
Suppose we have, besides the historical  prices, also the   High and Low prices  for each asset
(as we do in our experimental data). 
 Let $\tau >0$ be a number of periodic (monthly, in our experiments) observations of prices.
Then we estimate the value $Bid(t)$ of the bid at time $t$, as the mean of the 
daily Lows in the epoch $[t-\tau,t)$; correspondingly the value $Ask(t)$ of the ask at time $t$ is taken as the mean of the daily Highs in the epoch $[t-\tau,t)$.

We denote by $\bS(0)_b$ and $\bS(0)_{unb}$ the initial conservative market prices, which should lie in between these bid-ask intervals (i.e. in $[Bid(0), Ask(0)]$), using respectively the MEM Bounded 
(Eqs. \ref{dualent} and \ref{sol1})   
and the MEM Unbounded  method 
(Eqs. \ref{dualent2} and \ref{sol2}).

Table \ref{tb:bidask} shows these different estimates, at a given date, for the bid-ask intervals and the initial conservative market  price 
$\bS(0)$, for each of the stocks considered in our experiment. 
One can see that in practice both MEM methods (unbounded or bounded with $L=10^2$) give similar values for $\bS(0)$, and in both cases these values verify to be in, or close to, the interior of the Bid-Ask interval. Therefore, for our portfolio optimizations we shall only make use of the MEM unbounded method   for obtaining initial price 
$\bS^d(0)$.

\begin{table}[htp]
\begin{center}\begin{tabular}{|c|c|c|c|c|c|c|}\hline
Asset & $\widehat{\bS}(0)$ &   
 $\bS(0)_{b}$ ($10^2$) &   
 $\bS(0)_{unb}$ & $Bid$& $Ask$ \\ \hline
ALB & 20.04  & 
27.296 & 27.301 & 26.976 & 27.659 \\
AAPL & 13.53 &  
11.765 &  11.765 &  11.405 & 11.765 \\
APC & 40.07 &     
45.125 & 45.133 &  44.594 &  45.858 \\
GE & 16.37  &  
33.628 & 33.711 & 33.558 & 34.074    \\
WFC & 28.37 &   
 31.389 & 31.389 & 31.389 & 32.087\\
HAL & 17.38 &   
 29.563 & 29.569 & 29.168 & 29.968 \\
SHW &  59.91 &   
50.979 & 51.013 & 50.863 &   51.997 \\ \hline
\end{tabular} \caption{\small Price $\widehat{\bS}(0)$ in 2008-12-29,  
MEM-Bounded price $\bS(0)_{b}$ 
($L=10^2$), MEM-Unbounded price $\bS(0)_{unb}$, and bid-ask ranges.\label{tb:bidask}}
\end{center}
\end{table}

\subsubsection{Estimating $\bS(1)$}
The estimation $ \widehat{\bS}(1)$ of  $\bS(1)$, necessary  for solving the problem (\ref{utmin1}), can be done 
through a model of prices, for example 
a Geometric Brownian Motion (GBM) process,
using the information on estimated initial prices  $ \widehat{\bS}(0)$ and other parameters.
  We rather choose to estimate $\bS(1)$ by a data-driven 
 method, consisting on using a sufficiently long past history of gross returns up to time $t=0$ for each asset. Let ${\bf R} = (R_{ij})_{1\le i\le N,\, 1\le j\le M}$, where $R_{ij}$ is sample gross return of asset $j$ at time $i \le t=0$, 
and set
\[ \widehat{\bS}(1) = {\bf R}(\widehat{\bS}(0))^t\]
(where $t$ in the formula means transpose).

\subsubsection{Mean Variance Portfolio selection.}
For the sake of testing our portfolio methodology with other classical form of portfolio optimization, we consider optimizing portfolios with respect to the Markowitz's mean-variance utility function.
Recall that the Mean-Variance (long only) portfolio optimization is  given by the solution to the following problem:

\begin{equation}
\begin{aligned} 
& \underset{\bw \in \mathbb{R}^M}{\text{max}}
& & \mu ^t \bw - \frac{\gamma}{2}\bw^t\Sigma \bw\\
& \text{subject to} 
& &\langle \bw, \b 1\rangle = 1 , \;\; \ \bw \ge 0.
\end{aligned}
\end{equation}
where $\mu_i = E[X_i]$, $i=1,\ldots,M$, $\Sigma = \text{Var}(\bX)$ and $\gamma \ge 0$ is the risk aversion parameter. 
We observe that in most general cases of real data, $\mu$ is of the order of $10^{-2}$, and consequently $\Sigma \sim 10^{-4}$. Therefore, in order to built 
the Minimum Variance (MinVar) portfolio set  
$\gamma > 10^2$.

\subsubsection{General results}
To simulate portfolio's monthly returns we use 5 years of monthly data (60 observations), taken from the  dataset described above  with dates ranging  from 2000-01-01 to 2016-12-29. In order to obtain up to 13 different  scenarios (monthly periods where to assess the portfolio performance) we apply a rolling window analysis, each window period of 5 years, and using increments of 1 year of data between successive windows.
Figures \ref{fig:bp1} and \ref{fig:bp2}  show the portfolios  returns, for different monthly periods,  with respect to the method chosen to optimize the portfolio (exponential utility or mean-variance  maximization), and the way of   computing the initial values of portfolio 
(the $\bS(0)$).
We can see that for all but one case (the month ending in 2008-12-29, the date of the big market crash), optimizing the portfolios with initial prices computed by the MEM method improve returns substantially.
In Table \ref{tb:eumvmean} below we report the averages quartiles of monthly portfolio returns, over all 13 different months,  for portfolios optimized with respect to the exponential utility  (EU) objective and the minimum variance (MV) objective,   and taking $\bS(0)$ as current value ($\bS(0)_{current}$), 
or estimated by the MEM Unbounded method ($\bS(0)_{mem}$).
The results show that, in general, estimating initial prices $\bS(0)$ with the MEM method 
translates in better portfolio performance.

\begin{table}[!htp]
\begin{center}\begin{tabular}{|c|c|c|c|c|c|}\hline  \bf EU & $Q_0$ & $Q_1$ & Median &$Q_3$ & $Q_4$ \\\hline 
$\bS(0)_{current}$ & 0.923 & 0.985 & 1.015 & 1.045 & 1.118 \\\hline 
$\bS(0)_{mem}$ & 1.261 & 1.361 & 1.407 &1.451 & 1.561 \\\hline 
\end{tabular} 

\vspace{0.3cm}
\begin{tabular}{|c|c|c|c|c|c|}\hline  \bf MV & $Q_0$ & $Q_1$ & Median &$Q_3$ & $Q_4$ \\\hline 
$\bS(0)_{current}$ & 0.919 & 0.986 & 1.014 & 1.039 & 1.105 \\\hline 
$\bS(0)_{mem}$ & 1.251 & 1.352 & 1.392 &1.433 & 1.526 \\\hline 
\end{tabular} 
\caption{\small Average quartiles of monthly portfolio returns throughout all months considered\label{tb:eumvmean}}
\end{center}
\end{table}

We use the numerical software CVX with solver ECOS \citep{ecos}, and the  programming language R \citep{R}
to solve our optimization problems.

\begin{figure}[hbtp]
  \begin{center}
    \includegraphics[height=0.4\textwidth,width=0.49\textwidth]{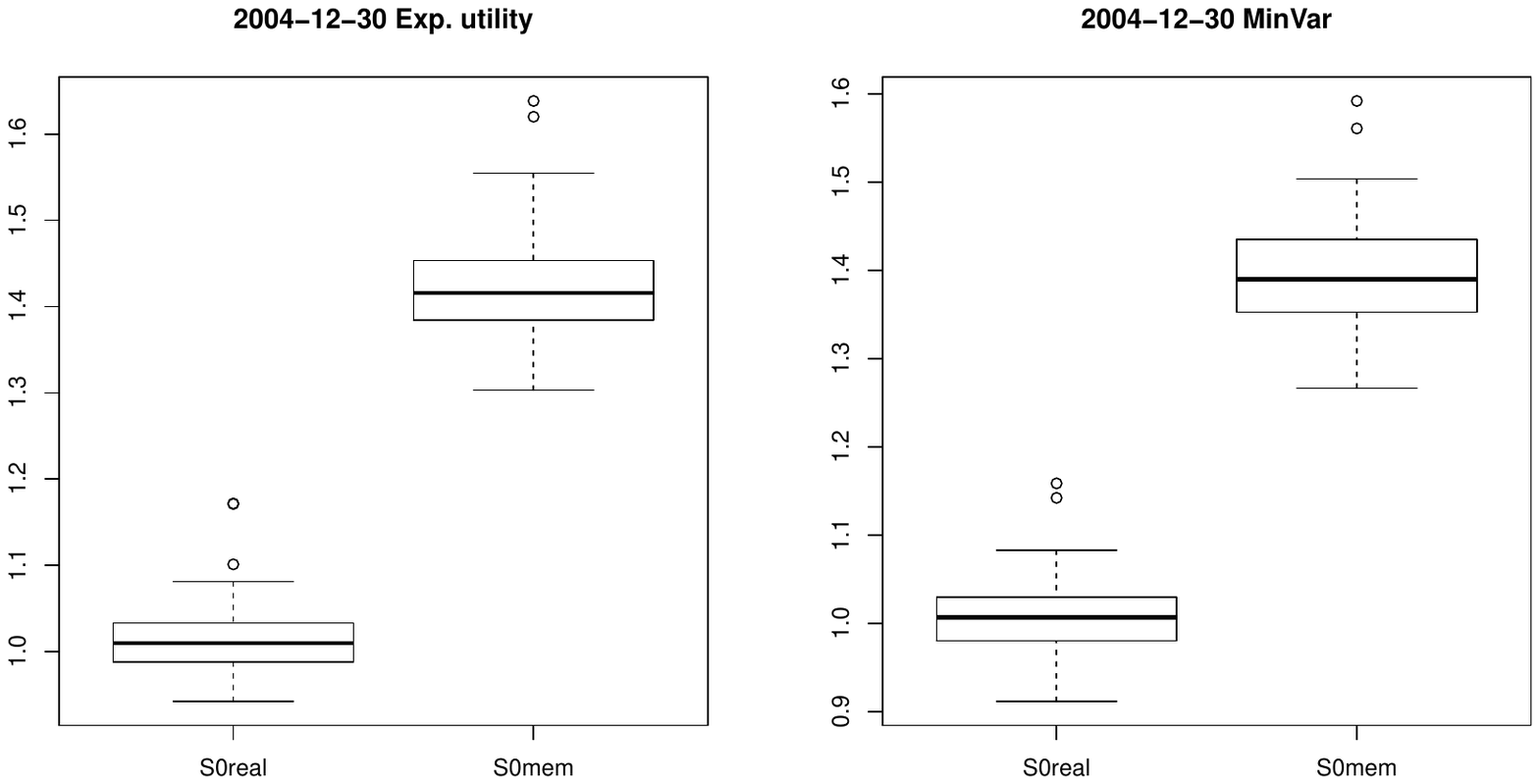}
     \includegraphics[height=0.4\textwidth,width=0.49\textwidth]{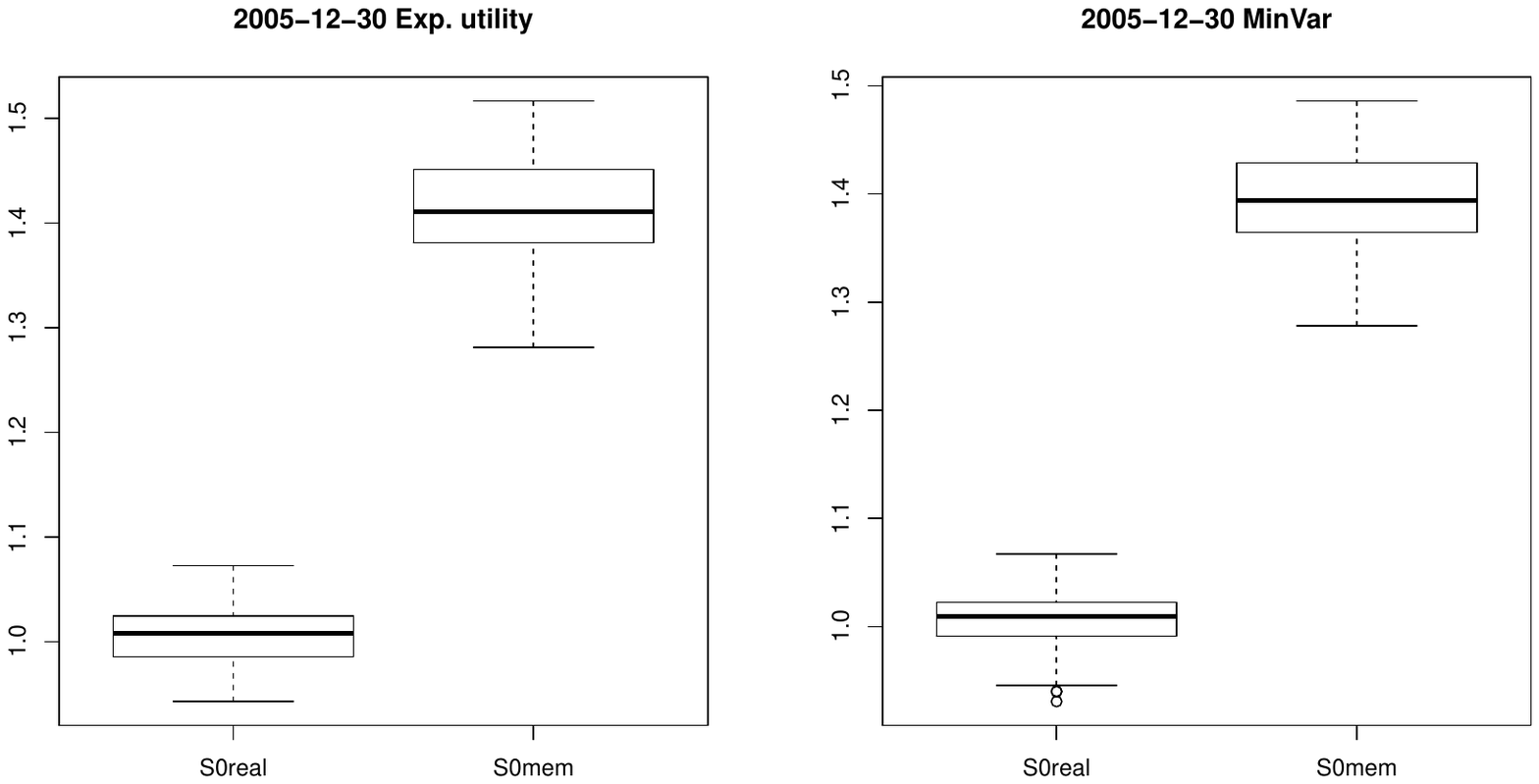}\\
       \includegraphics[height=0.4\textwidth,width=0.49\textwidth]{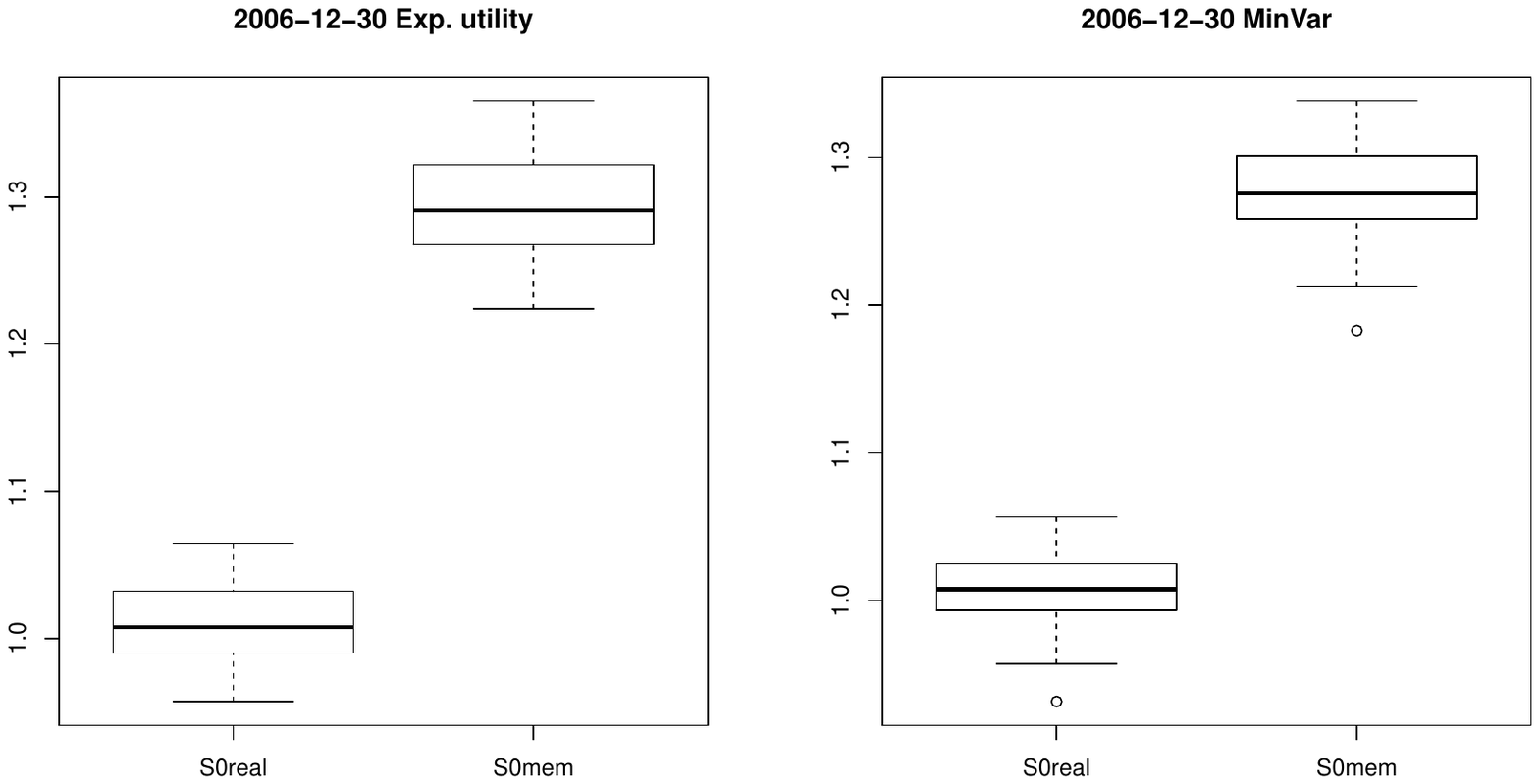}
     \includegraphics[height=0.4\textwidth,width=0.49\textwidth]{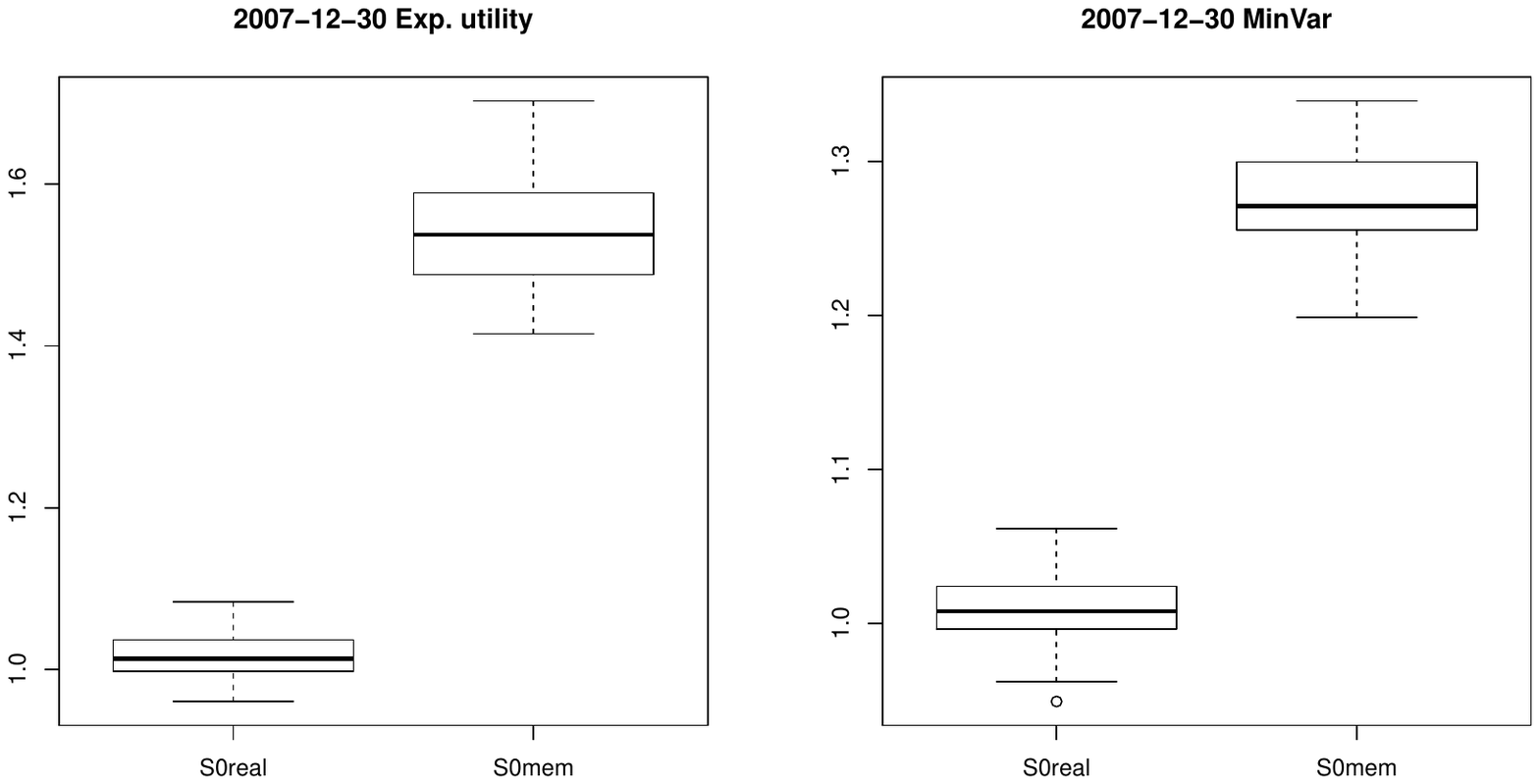}\\
      \includegraphics[height=0.4\textwidth,width=0.49\textwidth]{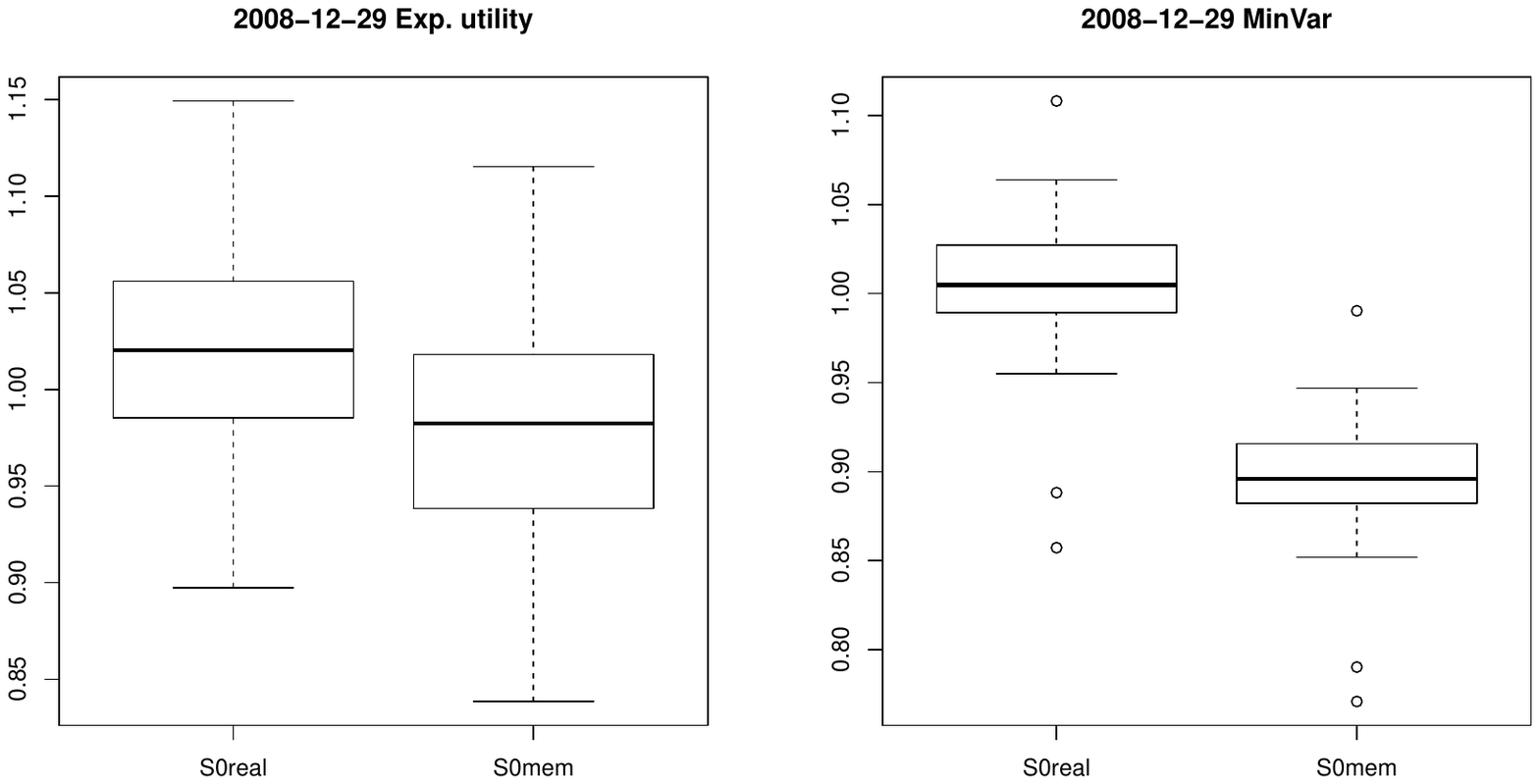}
     \includegraphics[height=0.4\textwidth,width=0.49\textwidth]{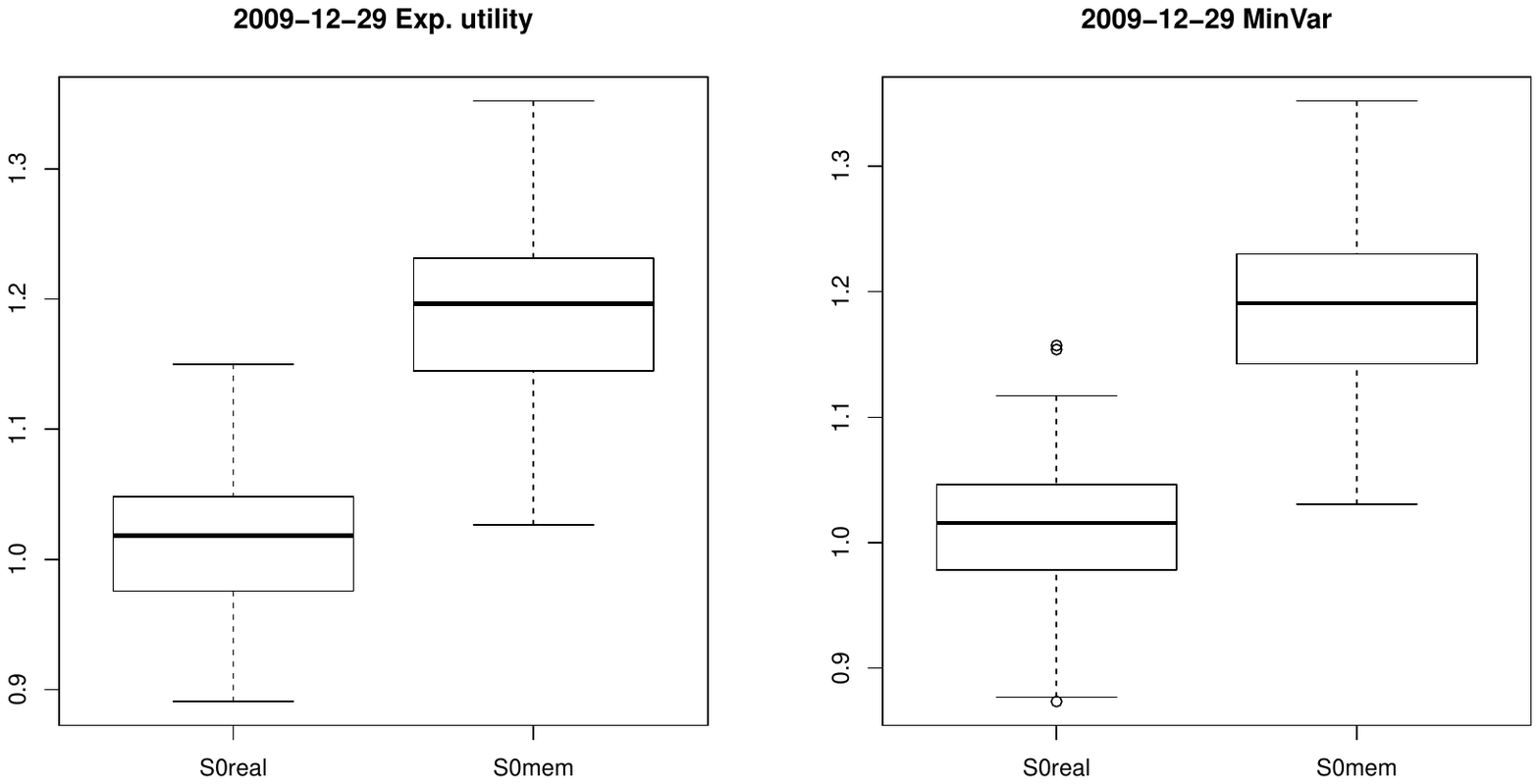}
       \caption{\footnotesize Boxplots of EU-type and MV-type portfolio's monthly returns (month end date indicated in header of each figure), for $\bS(0)$ obtained by: current price ($S0_{real}$), MEM method ($S0_{mem}$).\label{fig:bp1}}
  \end{center}
\end{figure}

\begin{figure}[hbtp]
  \begin{center}
    \includegraphics[height=0.4\textwidth,width=0.49\textwidth]{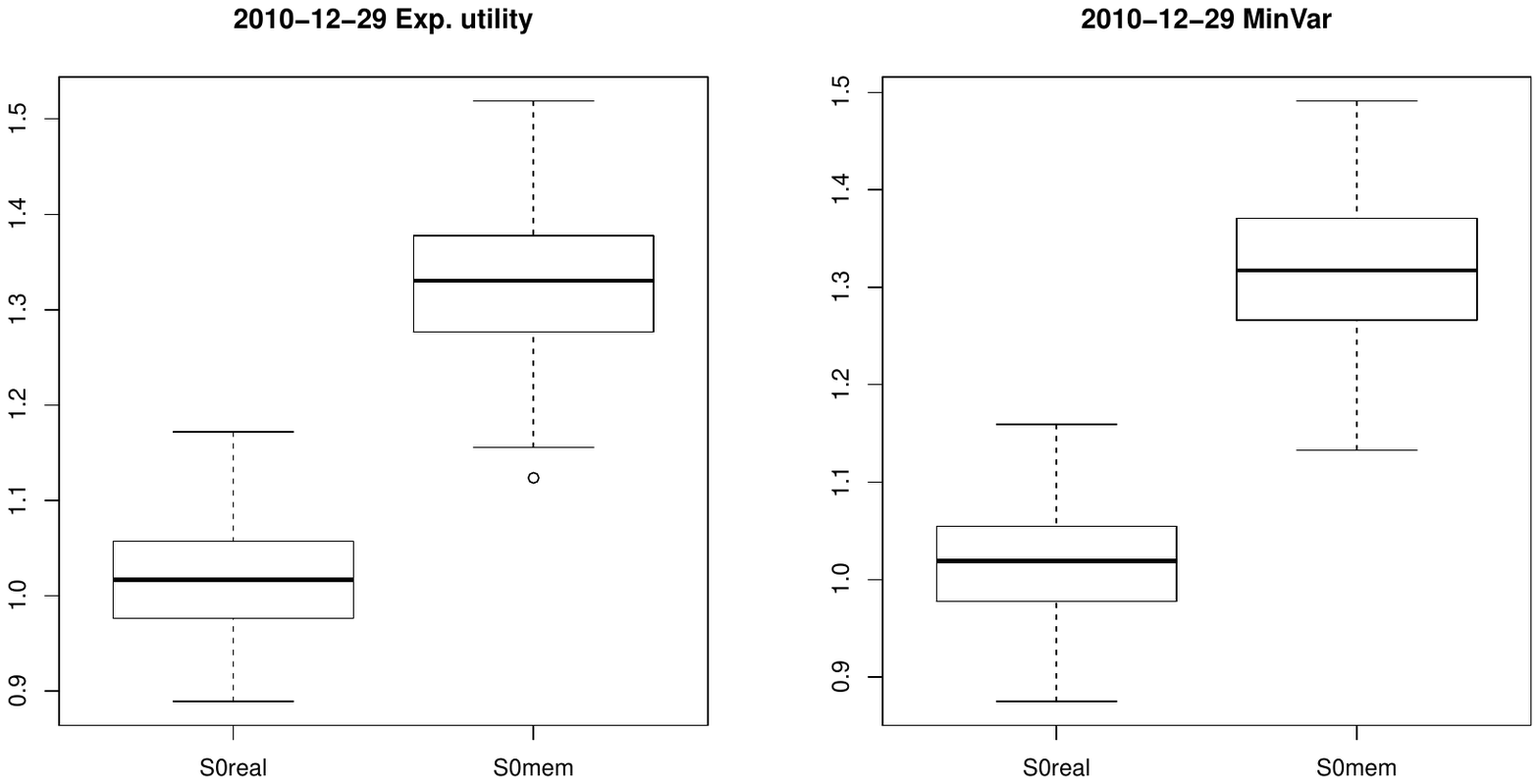}
     \includegraphics[height=0.4\textwidth,width=0.49\textwidth]{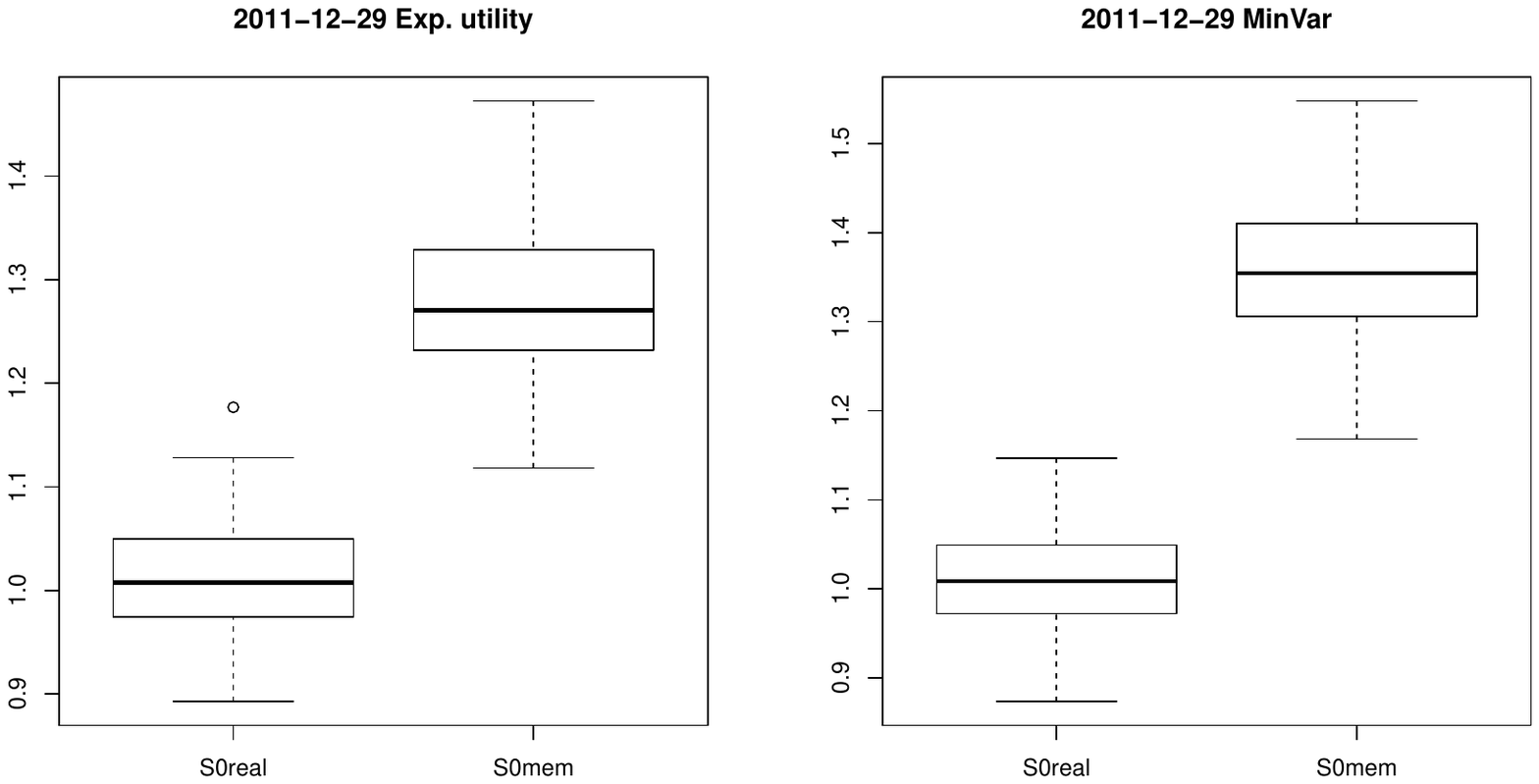}\\
       \includegraphics[height=0.4\textwidth,width=0.49\textwidth]{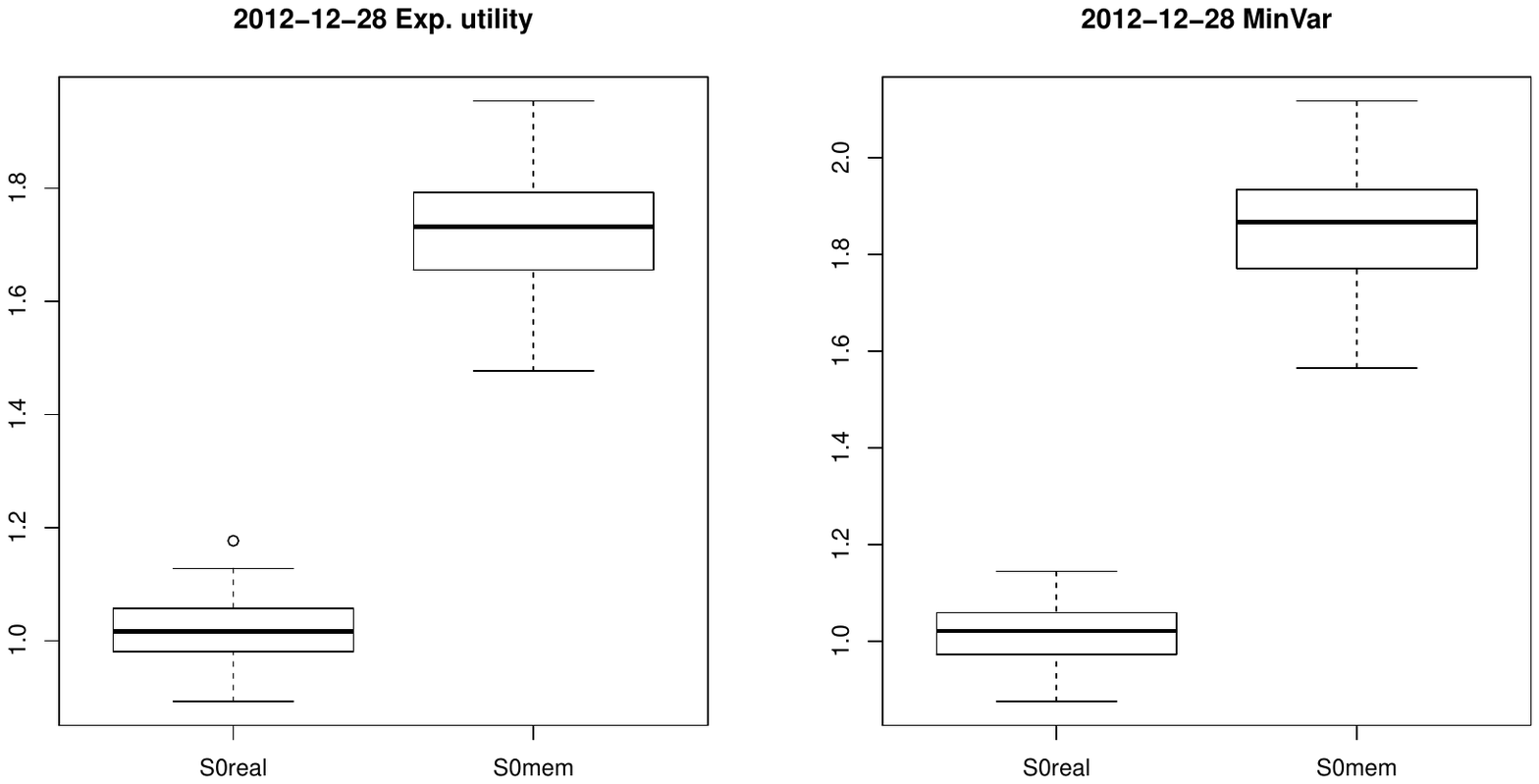}
     \includegraphics[height=0.4\textwidth,width=0.49\textwidth]{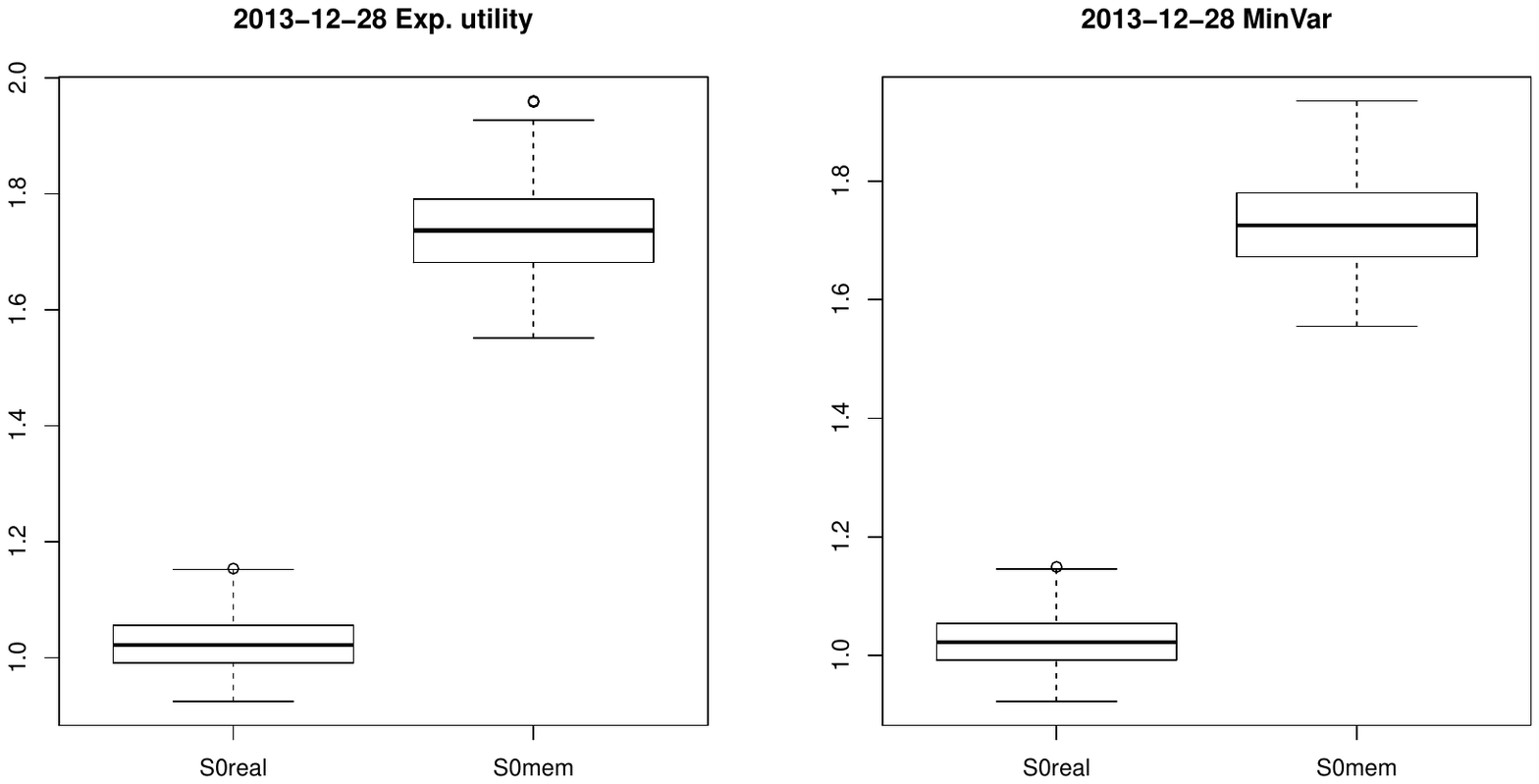}\\
      \includegraphics[height=0.4\textwidth,width=0.49\textwidth]{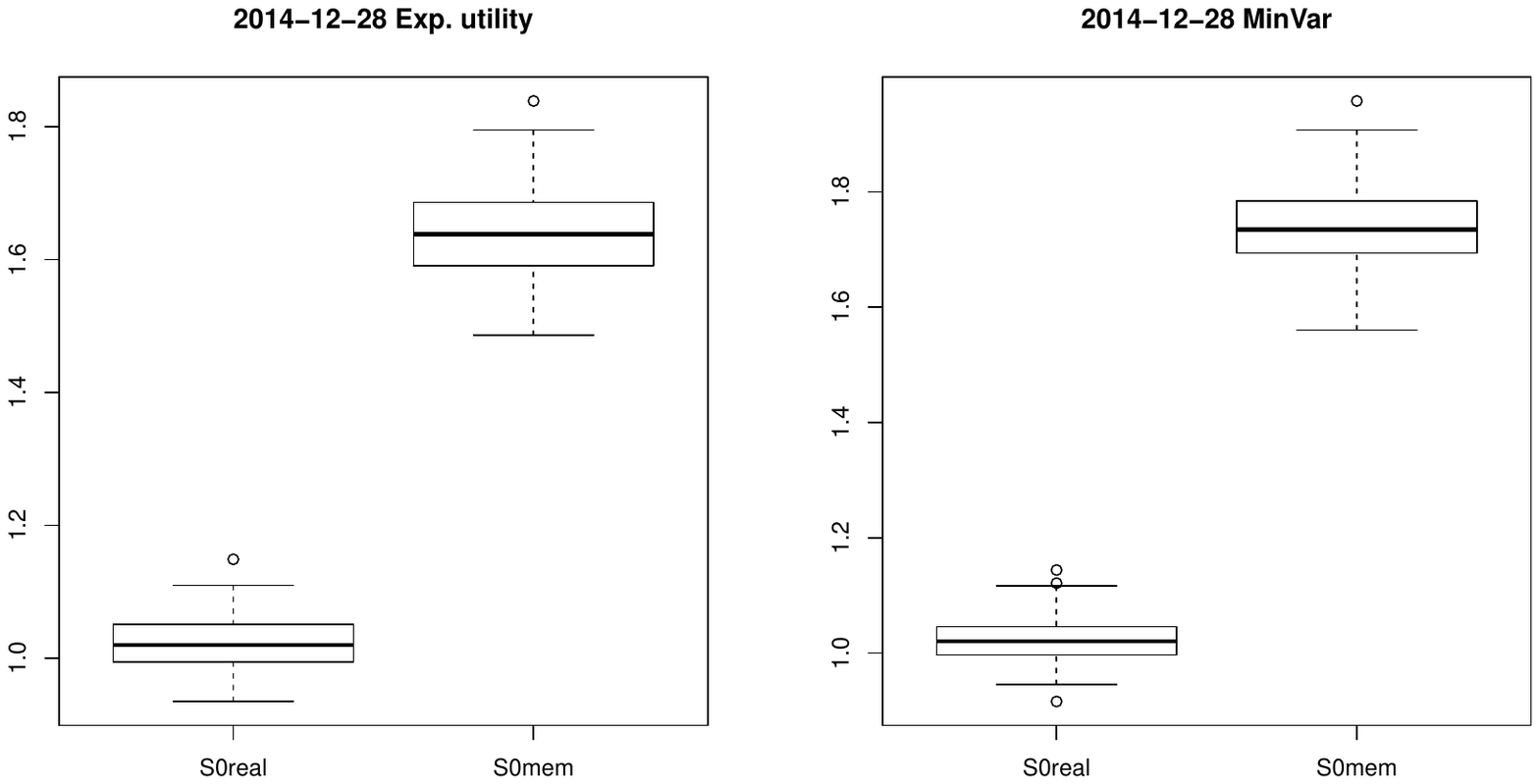}
     \includegraphics[height=0.4\textwidth,width=0.49\textwidth]{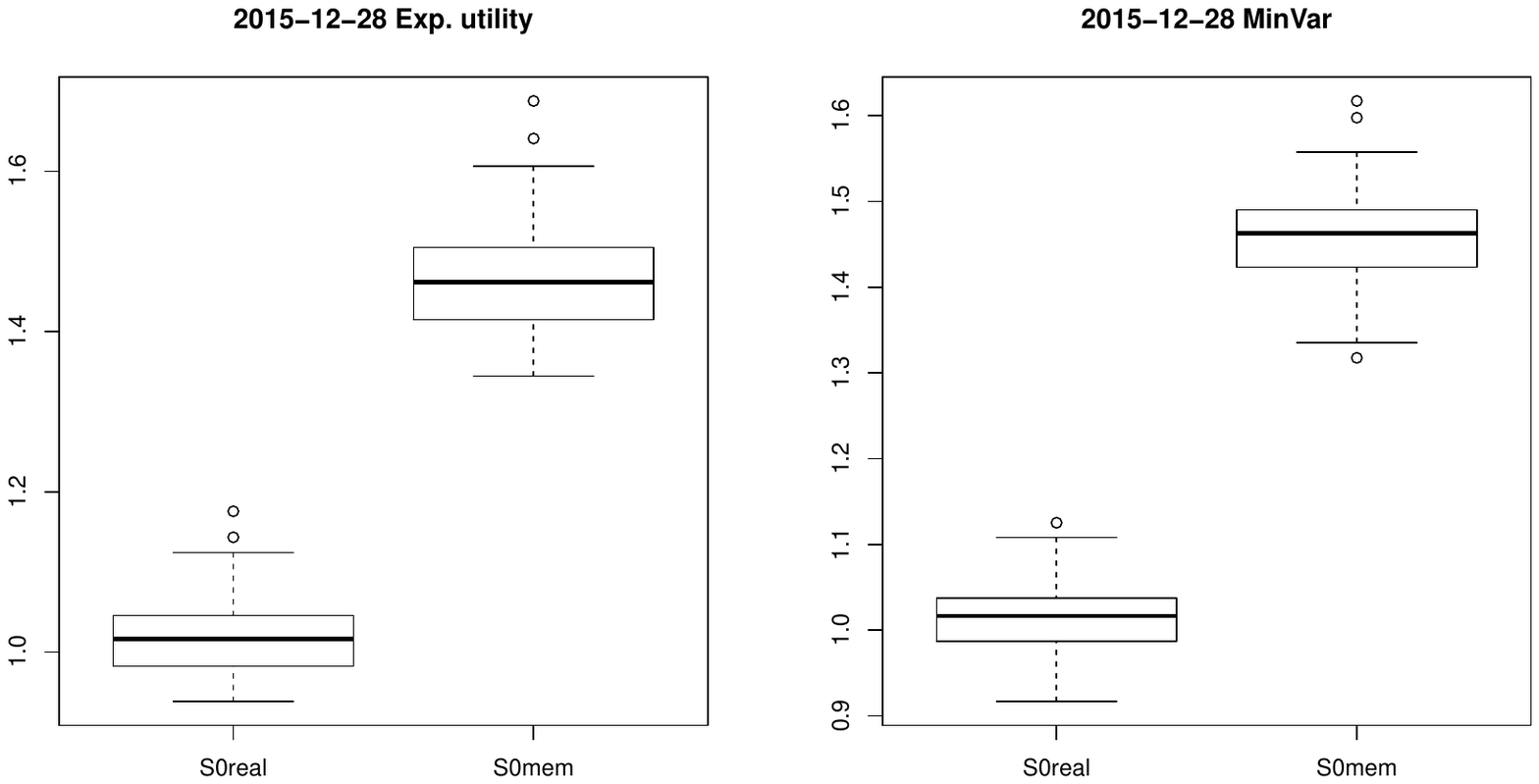}
       \caption{\footnotesize Continuation of Figure \ref{fig:bp1} showing comparison of  EU-type 
       and MV-type portfolio's monthly returns according to the way of computing  $\bS(0)$.\label{fig:bp2}}
  \end{center}
\end{figure}

 \section{Conclusions}\label{conc}
As said at the beginning, 
the underlying idea behind this paper is to provide a method to fill in a methodological detail in the process of portfolio optimization. Namely, that in order to comply with the cost constraint, we need to be able to assign prices to the   assets at $t=0,$ or in the notation used throughout, to be able to estimate $\bS(0)$, of which we formally know that
$\bS(0) = e^{-r}E_Q[\bS(1)],$  for $\bS(1)$   the random prices at $t=1$ and $Q$ some risk-neutral measure. 
We proposed an idea coming from insurance risk pricing, which turns out to be equivalent to the notion of risk neutral probability.

That part of the project produces nice results: that is, we are able to use a model free, non parametric approach to consistently estimate $\bS(0).$ 
The second part of the project, consists of using that estimate, or the constraint based on it, as part of an optimization process. 
Here a
 standard issue appears: How to model the tradeoff between risk and return, which is the initial starting point of every portfolio optimization problem.

We chose to exemplify the optimal portfolio configuration using two different methods: On one hand, a utility maximization procedure, and on the other a simple return-risk tradeoff in the Markowitz (volatility, return)-plane.
In both cases, it is clear the improvement in performance by  the portfolios   designed with  the initial prices of its constituent assets, estimated by the method based on maximum entropy in the mean that we have proposed in this paper.
 Finally, note that  the MEM method provides us with a non-parametric model free procedure to go around the lack of knowledge about the risk neutral probabilities.

\section{Appendices }\label{apps}
\subsection{Details about MEM}\label{A:sec1}
In this section we establish the basic formalism of the method of Maximum Entropy in the Mean (MEM). To begin, let us consider (\ref{prob3}) in a context free notation. The problem that we want to solve is to find $\bxi\in \mathcal{C}$  such that
$$ \bA\bxi \in \mathcal{K}$$
 where $\bA$ is an $M\times N$-matrix. 
We suppose that the sets $\mathcal{C}$ and $\mathcal{K}$ are convex and have non-empty interiors. Eventually we shall consider $\mathcal{C}=[0,L]^N$ or $[0,\infty)^N,$ and $\mathcal{K}=\prod_i[b_i,a_i].$  The parameter $L$ in our situation is related to a bound on the distortion function: it is given by the assumption that the values $\phi_j$, for $j=1, \ldots, N$, are bounded. If we do not want to preassign that value, we can set $L=\infty.$
 We shall give solutions for both cases of $L$, bounded or unbounded. 
 However, in our applications to portfolio optimization we shall see that  both cases of $L$ give similar results; 
 hence by the principle of parsimony one should prefer, in practice, to solve problem (\ref{prob3}) over  
 $\mathcal{C}=[0,\infty)^N$.
 
The idea behind MEM can be stated as follows: Consider an auxiliary random variable $\bX,$ taking values in $\mathcal{C}$ such that $\bxi = E_\mathbb{P}[\bX]$ with respect to an unknown probability measure $\mathbb{P}$ that is supported by $\mathcal{C}.$ If such probability measure can be found, the constraint upon $\bxi$ is automatically satisfied. To define $\bX$ we consider the probability space $(\mathcal{C}, \mathcal{B}(\mathcal{C}), \mathbb{Q}),$ where $\mathcal{B}(\mathcal{C})$ denotes the Borel subsets of $\mathcal{C},$ and $\mathbb{Q}$ is any ($\sigma$--finite) measure such that the convex hull $con(supp(\mathbb{Q}))$ generated by its support equals $\mathcal{C}.$ What this last requirement achieves for us, is that any probability $\mathbb{Q}$ having a strictly positive density $\rho(\bxi)$ with respect to $Q,$ will satisfy $E_\mathbb{P}[\bX]\in\mathcal{C}.$ The standard maximum entropy (SME) procedure enters as a stepping stone at this stage. On the class of densities (with respect to $\mathbb{Q}$) we define the entropy function by
$$S(\rho) = - \int_{\mathcal{K}}\rho(\bxi)\ln\rho(\bxi)d\mathbb{Q}(\bxi),$$
\noindent or $-\infty$ if the integral of $|\ln\rho(\bxi)|$ is not convergent. The choice of sign is conventional, we want to maximize entropy. The entropy function is strictly concave due to the strict concavity of the logarithm. To find a $\rho^*(\bxi)$ such that $d\mathbb{P}^*=\rho d\mathbb{Q}$ satisfies (\ref{prob1}) we solve the problem
\begin{equation}\label{prob12}
\mbox{Find}\;\;\;\rho^*\;\;\;\mbox{at which}\;\;\;\sup\{S(\rho)|\bA E_\mathbb{Q}[\bX]\in \mathcal{K}\} \;\;\;\mbox{is achieved}
\end{equation}
Since the maximization can be split into a sequence of two steps:
$$\sup_{\by\in\mathcal{K}}\sup\{S(\rho_{\by})|\bA E_P[\bX]=\by\}$$
\noindent in which the inner maximization has a standard solution  given by
\begin{equation}\label{rep1}
\rho_{\by}(\bxi) = \frac{e^{-\langle\blambda^*,\bA_e\bxi\rangle}}{Z(\blambda^*)}
\end{equation}
\noindent where the normalization factor is clearly given by
$$Z(\blambda) = \int_{\mathcal{C}}e^{-\langle\blambda^*,\bA\bxi\rangle}d\mathbb{Q}(\bxi).$$
The $\blambda^*$ appearing in (\ref{rep1}) is to be determined minimizing the convex function (dual entropy) given by
$$\Sigma(\blambda,\by) = \ln Z(\blambda) + \langle \blambda,\by\rangle$$
over $\{\blambda \in \mathbb{R}^M\,|\,Z(\blambda) < \infty\},$ which in all cases of interest for us will be $\mathbb{R}^M$ itself. Furthermore, and very important,
$$S(\rho_{\by}^*) = \Sigma(\blambda^*,\by).$$

Therefore, the double maximization described above becomes
\begin{equation}\label{minim}
\sup_{\by\in\mathcal{K}}S(\rho_{\by}^*) = \sup_{\by\in\mathcal{K}}\inf\{\Sigma(\lambda,\by)\,|\,\lambda\in\mathbb{R}^M\} = \inf_{\lambda\in\mathbb{R}^M}\{ \ln Z(\blambda) +\sup_{\by\in\mathcal{K}}\{\langle\blambda,\by\rangle\}\}.
\end{equation}

To make the pending computations explicit, note that
$$\sup_{\by\in\mathcal{K}}\{\langle\blambda,\by\rangle\} = \sum_{i=1}^M\left(\frac{a_i-b_i}{2}|\lambda_i| + \frac{a_i+b_i}{2}\lambda_i\right).$$

Now we shall complete the details in the two cases of interest to us, corresponding to $\mathcal{C}=[0,L]^N$ or $\mathcal{C}=[0,\infty)^N.$ In each case, we have to specify $Q$ and compute $Z(\blambda).$ 

\subsubsection{MEM: Bounded case}
For $L < \infty$, consider the following measure $\mathbb{Q}$ on $\mathcal{C}=[0,L]^N$:
$$\mathbb{Q}(d\bxi) = \prod_{j=1}^N\left(\epsilon_0(dx_j) + \epsilon_L(dx_j)\right).$$
We use the notation $\epsilon_a(dx)$ to denote the measure that puts a unit mass at the point $a.$  It takes a simple computation to verify that
$$Z(\blambda) = \prod_{j=1}^N\left(1 + e^{-(L\bA^t\blambda)_j}\right)$$
\noindent where the superscript ``t'' denotes transposition and $(\bA^t\blambda)_j = \sum_{i=1}^M\lambda_iA_{i,j}.$
With all this, the dual entropy to be minimized in (\ref{minim}) becomes
\begin{equation}\label{dualent}
\Sigma(\blambda) = \ln Z(\blambda) + \sum_{i=1}^M\left(\frac{a_i-b_i}{2}|\lambda_i| + \frac{a_i+b_i}{2}\lambda_i\right)
\end{equation}
Once the $\blambda^*$ that minimizes the right hand side of (\ref{dualent}) has been found, then
\begin{equation}\label{sol1}
\xi^*_j = L\frac{e^{-(L\bA^t\blambda^*)_j}}{1 + e^{-(L\bA^t\blambda^*)_j}}\;\;\;j=1, \ldots,N\\\\
\end{equation}
 
\medskip
At this point we add that when $a_i=y_i=b_i$ for $1=1, \ldots,M,$ that is when the data is pointwise, then the last term in (\ref{dualent}) becomes $\langle\blambda,\by\rangle.$ In this case the first order condition for $\blambda^*$ to be a minimizer of (\ref{dualent}) is equivalent to the fact that $\bA\bxi^* =\by$ as can be easily verified.  

To conclude the section we mention that the dual entropy $\Sigma(\blambda,\by)$ may not have a minimum. For example, $\Sigma(\blambda,\by)$ is bounded below in $\blambda$ for every $\by,$ but the infimum is reached at infinity; or when the problem is not well specified and data does not fall in the range of $\bA$, in which case $\Sigma(\blambda,\by)$ is not bounded below.  
For further details on this, see   \cite{BL} or   \cite{GV}.

\subsubsection{MEM: Unbounded case} 
The problem that we want to solve now is to find $\bxi\in \mathcal{C}=[0, \infty)^N$  such that
$$ \bA\bxi \in \mathcal{K}$$
\noindent where $\bA$ is an $M\times N$-matrix. As in the bounded case, $\mathcal{C}$ and $\mathcal{K}$ are convex and have non-empty interiors. Everything goes almost verbatim as above, except that the reference measure that we consider now is a product of un-normalized Poisson measures:
$$d\mathbb{Q} = \prod_{j=1}^N \left(\sum_{n=1}^\infty \frac{1}{n!}\epsilon_n\; d\xi_j \right).$$
This time the function $Z(\blambda)$ is given by
$$Z(\blambda) = \prod_{j=1}^N \exp\Big(e^{-(\bA^t\blambda)_j}\Big)$$
To continue as above we form the function $\Sigma(\blambda)$ and obtain
\begin{equation}\label{dualent2}
\Sigma(\blambda) = \sum_{j=1}^N e^{-(\bA^t\blambda)_j} + \sum_{i=1}^M\left(\frac{a_i-b_i}{2}|\lambda_i| + \frac{a_i+b_i}{2}\lambda_i\right)
\end{equation}
Once again, when the $\blambda^*$ that minimizes the right hand side of (\ref{dualent2}) has been found, then the solution to $ \bA\bxi \in \mathcal{K}$ is given by 
\begin{equation}\label{sol2}
\xi^*_j = e^{-(\bA^t\blambda^*)_j}, \;\;\;j=1, \ldots,N\\\\
\end{equation}
Notice as well that when $a_i=y_i=b_i$ then
$$\Sigma(\blambda) = \sum_{j=1}^N e^{-(\bA^t\blambda)_j} + \sum_{i=1}^M\lambda_i y_i.$$
In this case, the first order condition on $\blambda^*$ to be a minimizer looks like
$$-\sum_{j=1}^N A_{i.j}e^{-(\bA^t\blambda^*)_j} + y_i\;\;\;\;i=1,...,M$$
from which it is clear that(\ref{sol2}) satisfies the equation $\bA\bxi = \by.$

\subsubsection{Sub-differential of (\ref{dualent}) }
As a first step towards a solution, we have the following:

\begin{lemma}
With the notations introduced above, the function $\Sigma(\blambda)$ defined in (\ref{dualent}) 
is strictly convex
\end{lemma}

Observe that 
$\partial_i\Sigma(\blambda) := \partial\Sigma(\blambda)/\partial\lambda_i$ is defined except at $\lambda_i=0$, where it is only sub-differentiable \citep{BL}. More explicitly, 

\begin{equation}\label{gradS1}
\partial_i\Sigma(\blambda)= \left\{ \begin{array}{l@{\quad}l}
\partial \ln Z(\blambda)/\partial\lambda_i + a_i\ & \mbox{when }\ \lambda_i >0, \\ \\
 \partial \ln Z(\blambda)/\partial\lambda_i + b_i\ & \mbox{when }\ \lambda_i <0, \\ \\
 \partial \ln Z(\blambda)/\partial\lambda_i \in (b_i, a_i)\ & \mbox{when }\ \lambda_i = 0,
\end{array} \right.
\end{equation}

Moreover,  in the bounded case,
$\ln Z(\blambda) = \sum_{j=1}^N\ln\left(1 + e^{-(L\bA^t\blambda)_j}\right)$, and using the expression
$(\bA^t\blambda)_j = \sum_{k=1}^M\lambda_kA_{k,j}$, we have that

\begin{equation}\label{gradS2}
\partial \ln Z(\blambda)/\partial\lambda_i  = -\sum_{j=1}^N LA_{i,j}\frac{ e^{-(L\bA^t\blambda)_j}}{1+ e^{-(L\bA^t\blambda)_j}}
\end{equation}
and in the unbounded case
$$\ln Z(\blambda) = \sum_{j=1}^N e^{-(\bA^t\blambda)_j}$$
and we  have 
\begin{equation}\label{gradS3}
\partial \ln Z(\blambda)/\partial\lambda_i  = -\sum_{j=1}^N A_{i,j}e^{-(\bA^t\blambda)_j}
\end{equation}
Finally, we have the general statement
\begin{theorem}
Suppose that the infimum $\blambda^*$ of $\Sigma(\blambda)$ is reached at some $\blambda^*$ in the interior of $\{\blambda \in \mathbb{R}^M : Z(\blambda) < \infty\}$. Then 
$$\partial_i \Sigma(\blambda) \in \{0\}, \quad i=1, \ldots, M.$$ 
\end{theorem}

\subsection{Distortion functions and risk neutral densities: General case}\label{A:sec2}
Here we shall extend the result presented in     Section \ref{DF_RNP} to the multidimensional case. Remember that we are considering an $M$-dimensional, positive random vector $\bS(1)$ describing the values of the assets at $t=1.$ Recall as well that we are supposing each $S_i(1)$ to be continuously distributed having cumulative distribution function $F_i(x)$ with strictly positive densities.  Let us denote by $F(\bx)=P(S_1(1)\leq x_1, \ldots, S_M(1)\leq x_M)$ the joint distribution function of $\bS(1).$ Let $g$ be the distortion function introduced in     Section \ref{DF_RNP} and again $\psi(u)=1-g(1-u).$ We begin by establishing some properties of the  distorted cumulative distribution function $G(\bx)=\psi(F(\bx))$.
 Let us denote by $G_i(x_i)=\psi(F_i(x_i))$ the one dimensional distorted cumulative density of $S_i(1)$.
\begin{lemma}\label{abscont}
With the notations introduced above, $G(\bx)$ has marginals $G_i(x_i).$
\end{lemma}
The proof for the general case is similar to that of the 2-dimensional case, but notationally simpler.
Notice, for example, that, for fixed $x_1\in[0,\infty)$
$$\int_0^\infty d_2G(x_1,x_2)=\int_0^\infty d_2\psi(F(x_1,x_2))$$
$$=\psi(F(x_1,\infty))-\psi(F(x_1,0))=\psi(F_1(x_1))=G_1(x_i).$$
From this it is clear that, for any positive or bounded, measurable $h$ defined on $[0,\infty),$ we have $E[h(S_i)]=E^\psi[h(S_i)],$ and in particular,
\begin{equation}\label{RN2}
e^{-r}E^\psi[S_i] = e^{-r}\int_{\mathbb{R}_{+}^M} x_i d^MG(\bx) = e^{-r}\int_0^\infty x_i d_i G_i(x_i),
\end{equation}
where $d^M := d_1d_2\cdots d_M$. 
Therefore, we could have stated Problem  (\ref{prob0}) as consisting of finding a distortion function $g$ such that
\begin{equation}\label{prob4}
e^{-r}E^{\psi}[S_i] = e^{-r}\int_{\mathbb{R}_{+}^M} x_id^M\psi(F(\bx)) \in (b_i,a_i),\;\;\;\mbox{for}\;\;i=1, \ldots,M.
\end{equation}
To continue, according to Sklar's theorem \citep{Sklar59}, there exists  a (unique under our assumptions) $N-$dimensional distribution function $C_F(\bu)$ on $[0,1]^M$ such that $F(\bx) = C_F(F_1(x_1),$ \ldots, $F_M(x_M)).$ 
We shall now verify that the co\-pu\-la $C_G$ of the distorted distribution function is a $\psi$-distorted 
co\-pu\-la of $C_F$. See   \cite{durante} for an introduction to the subject.
\begin{proposition}\label{distcop}
 With the notations introduced in the previous paragraph we have,
$$C_G(\bu) = \psi\big(C_F(\psi^{-1}(u_1), \ldots,\psi^{-1}(u_M))\big).$$
\end{proposition}
The proof goes as follows: Recall that $G_i(x_i)=\psi(F_i(x_i)),$ which implies that 
$F_i(x_i)=\psi^{-1}(G_i(x_i)).$
Therefore
\begin{eqnarray*}
G(\bx) &=& \psi\big(C_F(F_1(x_1), \ldots,F_M(x_M))\big)\\
 &=& \psi\big(C_F(\psi^{-1}(G_1(x_1)), \ldots,\psi^{-1}(G_M(x_M)))\big)
\end{eqnarray*}
Since $G(\bx)= C_G\big(G_1(x_1), \ldots,G_M(x_M)\big)$, the desired result follows by setting $G_i(x_i)=u_i.$\\
In order to relate $G(\bx)$ to a risk neutral density, we have to establish the absolute continuity of $G(\bx)$ with respect to $\prod_{i=1}^Md_iF_i(x_i).$ To motivate the result, consider the following computations.\\
For $M=1,$
$$d_1G(x_1) = \psi^{(1)}(F_1)d_1F_1(x_1).$$
For $M=2,$ drop the subscript $F$ from $C_F$ and $C_i$ denotes $\partial C/\partial u_i,$ etc. 
The chain rule applied to $G(\bx) = \psi(C(F(\bx))$ yields
\begin{eqnarray*}
d_2d_1G(x_1,x_2) &=& \Big(\psi^{(2)}(C(F))C_1(F)C_2(F) \\
  & &\ +\ \psi^{(1)}(C(F))C_{1,2}(F)\Big)d_1F_1(x_1)d_2F_2(x_2).
\end{eqnarray*}
For $M=3,$
$$  d_3d_2d_1G(\bx) = \Big(\psi^{(3)}(C(F))C_1(F)C_2(F)C_3(F) + \psi^{(2)}\big(C_{1,2}(F)C_3(F)$$
 $$+\ C_{2,3}(F)C_1(F)+C_{1,3}(F)C_2(F)\big) + C_{1,2,3}(F)\Big)d_1F_1(x_1)d_2F_2(x_2)d_3F_3(x_3).$$
Note that the symmetry with respect to the exchange of the subscripts is related to the fact that we can integrate with respect to any variable in any order.
We shall now split Theorem \ref{RNP-DF} into two parts (to keep its length under control). But let us first introduce some more notations. Let $[M]=\{1, \ldots,M\}$ and put
$$\mathcal{P}(M,k) = \{\mbox{Partitions of}\;\;\;[M]\;\;\;\mbox{into}\;\;\;k\;\;\mbox{blocks}\}.$$
For each subset $J=\{i_1,i_2, \ldots,i_n\}$ of $[M]$ we write
$$C_J = \frac{\partial^N C}{\partial u_{i_1} \cdots \partial u_{i_n}}.$$

\begin{theorem}\label{RNM-DF2}
With the notations introduced above, suppose that $\psi(u)$ is $M$-times continuously differentiable with positive derivatives $\psi^{(k)}.$ Suppose that $C$ has continuous derivatives $C_J$ for all $J\subset [M].$ (These are non-negative since $C$ is a cumulative distribution function and we differentiate with respect to different variables). Then we have
\begin{equation}\label{df1}
d^M G(\bx) = \Big[\sum_{k=1}^M  \psi^{(k)}(C(F))\Big(\sum_{J\in\mathcal{P}(M,k)}C_J(F)\Big)\Big]\prod_{i=1}^Md_iF_i(x_i).
\end{equation}
That is, $G(\bx)$ is absolutely continuous with respect to $\prod_{i=1}^Md_iF_i(x_i),$ with density given by the term in square brackets. Denote that by $\rho(F(\bx)).$  Furthermore 
(consider $\psi(u)=u$ if you want),
$d^M F(\bx) = C_{[M]}(F)\prod_{i=1}^M d_iF_i(x_i).$
\end{theorem}
 
\begin{corollary}\label{coro1}
If both $\rho(F)>0$ and $C_{[M]}>0,$ then $G \sim F$ with density $\frac{\rho(F(\bx))}{C_{[M]}(F(\bx))}.$
Furthermore, setting $\Lambda=\frac{\rho(F(\bS))}{C_{[M]}(F(\bS))}$ and $dQ=\Lambda dP,$ then for any measurable, positive or bounded $h(\bx)$, we have
$E_Q[h(\bS)] = E^{\psi}[h(\bS)]$, 
and in particular
$$e^{-r}E_Q[S_i] \in (b_i,a_i), \;\;\;\mbox{for}\;\;\;i=1, \ldots, M.$$
\end{corollary}

The statement of the converse will come at the end of the following observations. Note that if $Q$ is a probability on $(\Omega,\mathcal{F})$ which is equivalent to $P,$ if all functions that we are going to integrate are $\sigma(\bS)-$measurable, we might suppose that there is a function $\Lambda:\mathbb{R}_{+}^M\rightarrow (0,\infty)$ such that $d^MG(\bx)=\Lambda(\bx)d^MF(\bx).$  Here $G(\bx)=Q(\bS \leq \bx).$
We shall suppose that $d^MF(\bx)$ is equivalent to $\prod_{i=1}^M d_iF_i(x_i)$ with strictly positive density given by $C_M(F(\bx)).$
Write
 $\Lambda(\bx)=k(F_1(x_1), \ldots,F_M(x_M))$,
  with $k(u_1, \ldots,u_M)=\Lambda(F_1^{-1}(u_1), \ldots,F^{-1}_M(u_M)),$ then
$$d^MG(\bx) = k(F_1(x_1), \ldots,F_M(x_M))C_M(F_1(x_1), \ldots,F_M(x_M))\prod_{i=1}^M d_iF_i(x_i).$$
Note as well that the term in square brackets in (\ref{df1}) is obtained successively differentiating $\psi(C(\bu))$ with respect to $u_1, \ldots,u_M.$ Thus to obtain $\psi$ from $k(\bu)C_M(\bu)$ all we have to do is to consider $d_1 \cdots d_M\psi(C(\bu)) = k(\bu)C_M(\bu)$ as an equation to be solved for $\psi$ by successive integrations and set
$$\psi(t) = \int_0^t\left(\int_0^1\left(\int_0^1 \ldots\left.\left(\int_0^1k(\bu)C_M(\bu)du_M\right)du_{M-1}\right) \ldots\right)du_2\right)du_1.$$
When carrying out the integrations, one has to keep in mind that 
$$C^{(M-k)}(u_1, \ldots,u_{M-k})=C(u_1, \ldots,u_{M-k},1,1, \ldots,1)$$
 is an $M-k$-dimensional copula obtained by setting the last $k$ components equal to $1,$ and that $C^{(1)}(u_1)=u_1.$
We collect the previous  remarks under the following theorem. 
\begin{theorem}\label{RNM-DF3}
With the notations just introduced, suppose that there exists a probability $Q\sim P$ such that
$$e^{-r}E_Q[S_i] \in (b_i,a_i),\;\;\; \mbox{for}\;\;i=1, \ldots, M.$$ 
Let $G(\bx)=Q(\bS\leq\bx)$ and suppose that 
$$d^MG(\bx) = k(F_1(x_1), \ldots, F_M(x_M))C_M(F_1(x_1), \ldots,F_M(x_M))\prod_{i=1}^M d_iF_i(x_i).$$
Then there is an increasing function $\psi(t)$ defined on $[0,1]$ by
$$\psi(t) = \int_0^t\int_0^1\int_0^1 \ldots\int_0^1k(\bu)C_M(\bu)du_M\, du_{M-1}\, \ldots\, du_2\, du_1,$$
\noindent which necessarily satisfies $\psi(0)=0$ and $\psi(1)=1$ such that $G(\bx)=\psi(F(\bx)),$ and therefore
$$e^{-r}E^\psi[S_i] \in (b_i,a_i),\;\;\;\mbox{for}\;\;i=1, \ldots, M.$$ 
\end{theorem}

\begin{acknowledgements}
 The authors gratefully acknowledge the fruitful comments of the anonymous reviewers  that lead to an improvement of the original manuscript.

\noindent
A. Arratia research is supported by  grant TIN2017-89244-R from MINECO (Ministerio de Econom{\'\i}a, Industria y Competitividad) and the recognition 2017SGR-856 (MACDA) from AGAUR (Generalitat de Catalunya).  
\end{acknowledgements}

\end{document}